\documentclass[a4paper,english]{IEEEtran}
\usepackage[T1]{fontenc}
\usepackage{float}
\usepackage{multirow}
\usepackage{amsbsy}
\usepackage{amstext}
\usepackage{amsthm}
\usepackage{graphicx}
\usepackage{epstopdf}
\makeatletter

\pdfpageheight\paperheight
\pdfpagewidth\paperwidth

\providecommand{\tabularnewline}{\\}
\floatstyle{ruled}
\newfloat{algorithm}{tbp}{loa}
\providecommand{\algorithmname}{Algorithm}
\floatname{algorithm}{\protect\algorithmname}

\theoremstyle{plain}
\newtheorem{thm}{\protect\theoremname}
\theoremstyle{plain}
\newtheorem{prop}[thm]{\protect\propositionname}

\usepackage{amsmath}
\usepackage{amssymb}
\usepackage[numbers,sort]{natbib}

\usepackage{graphicx}
\usepackage{subfigure}
\usepackage{pdfsync}
\usepackage{url}
\usepackage{pdfsync}
\usepackage{todonotes}
\usepackage{tikz}
\def\checkmark{\tikz\fill[scale=0.4](0,.35) -- (.25,0) -- (1,.7) -- (.25,.15) -- cycle;}

\makeatother

\usepackage{babel}
\providecommand{\propositionname}{Proposition}
\providecommand{\theoremname}{Theorem}

\begin{document}

\title{Energy Efficiency Optimization in MIMO Interference Channels: A Successive
Pseudoconvex Approximation Approach}

\author{Yang Yang, Marius Pesavento, Symeon Chatzinotas and Björn Ottersten\thanks{Y. Yang, S. Chatzinotas and B. Ottersten are with Interdisciplinary
Centre for Security, Reliability and Trust, University of Luxembourg,
L-1855 Luxembourg (email: yang.yang@uni.lu, symeon.chatzinotas@uni.lu,
bjorn.ottersten@uni.lu). Their work is supported by FNR projects SATSENT
and CIPHY.}\thanks{M. Pesavento is with Communication Systems Group, Technische Universit\"{a}t
Darmstadt, 64283 Darmstadt, Germany (email: pesavento@nt.tu-darmstadt.de).
His work is supported by the EXPRESS Project within the DFG Priority
Program CoSIP (DFG-SPP 1798) }}
\maketitle
\begin{abstract}
In this paper, we consider the (global and sum) energy efficiency
optimization problem in downlink multi-input multi-output multi-cell
systems, where all users suffer from multi-user interference. This
is a challenging problem due to several reasons: 1) it is a nonconvex
fractional programming problem; 2) the transmission rate functions
are characterized by (complex-valued) transmit covariance matrices;
and 3) the processing-related power consumption may depend on the
transmission rate. We tackle this problem by the successive pseudoconvex
approximation approach, and we argue that pseudoconvex optimization
plays a fundamental role in designing novel iterative algorithms,
not only because every locally optimal point of a pseudoconvex optimization
problem is also globally optimal, but also because a descent direction
is easily obtained from every optimal point of a pseudoconvex optimization
problem. The proposed algorithms have the following advantages: 1)
fast convergence as the structure of the original optimization problem
is preserved as much as possible in the approximate problem solved
in each iteration, 2) easy implementation as each approximate problem
is suitable for parallel computation and its solution has a closed-form
expression, and 3) guaranteed convergence to a stationary point or
a Karush-Kuhn-Tucker point. The advantages of the proposed algorithm
are also illustrated numerically.
\end{abstract}


\section{Introduction}

In the era of 5G and Internet of Things by 2020, the number of connected
devices is predicted to reach 50 billions \cite{IMT}. On one hand,
as compared to current systems, the data rate should be 1000x higher
to serve these devices simultaneously. On the other hand, the significant
increase in the data rate is expected to be achieved at the same or
even a lower level of energy consumption. Therefore the so-called
energy efficiency (EE) is a key performance indicator that should
be considered in the design of transmission schemes.

In this paper, we adopt the notion of EE as the ratio between the
transmission rate and the consumed energy, which has a unit of bits/Joule,
and we study the EE maximization problem in a downlink multi-input
multi-output (MIMO) multi-cell system, where the base stations (BSs)
are transmitting in the same frequency band to allow full frequency
reuse and the users suffer from multi-user interference. This problem
is challenging due to several practical difficulties:
\begin{enumerate}
\item [(D1)]The transmission rate in the interference channel is a nonconcave
function of the transmit covariance matrices.
\item [(D2)]The energy consumption depends not only on the transmission
power but also on the processing power that increases with the transmission
rate.
\item [(D3)]In MIMO systems, the transmission rate functions are characterized
by (complex-valued) transmit covariance matrices.
\end{enumerate}
For the sake of an intuitive understanding of the challenging nature,
consider the sum rate maximization problem, which is a special case
of the EE optimization problem if the power consumption is a constant:
it has been proved in \cite{Luo2008} that the sum rate maximization
problem in interference channels is nonconvex and finding its globally
optimal point is NP-hard. Due to the high complexity of global optimization,
we are mainly interested in iterative algorithms with parallel implementations
that can efficiently find stationary points.

In a multi-cell network, multiple transmission links coexist that
negatively influence each other through the multi-user interference.
The conflicting interests of different links make the EE maximization
problem a multi-objective optimization problem and there are several
commonly adopted design metrics with different rationale. For example,
the global energy efficiency (GEE), which is defined as the ratio
between the sum transmission rate and the total power consumption,
is a meaningful measure for the EE of the whole network. Nevertheless,
it may not be relevant in a heterogeneous network, where different
transmission links may have different priorities. The EE of this network
is better captured by the (weighted) sum energy efficiency (SEE),
defined as the sum of all individual EE.

\emph{Related work.} The EE optimization problem has received considerable
attention in recent years and it has been studied from different perspectives.
For example, to address (D1), orthogonal transmission schemes based
on user selection or interference cancellation, are adopted in some
of the early works \cite{Ng2012,Ng2012a,Xu2013,Xu2014} so that the
transmission rate functions are concave in the transmit covariance
matrices. However, this scheme is not optimal due to the inefficient
reuse of spectrum, especially considering the large number of devices
in future networks and the existing frequency bandwidth limitations.

Along the direction of nonorthogonal multiple access, the GEE maximization
in MISO systems has been studied in \cite{Tervo2015}, where the authors
considered additional Quality-of-Service (QoS) constraints, in terms
of each link's guaranteed minimum transmission rate. The SEE optimization
problem with QoS constraints is studied in \cite{Tervo2017}. Compared
with the GEE function, the SEE function is more difficult to optimize
because it is the sum of multiple fractional functions, while each
individual fractional function is the ratio of a nonconcave function
and a nonconvex function. In MISO systems, the transmission rate is
a function of the SINR which is a scalar quantity and the algorithms
proposed in \cite{Nguyen2015,Tervo2015,Tervo2017} are built upon
this property. Thus they are not applicable for MIMO systems where
the transmission rate is a function of the (complex-valued) transmit
covariance matrices.

The sequential pricing algorithm for SEE maximization in MISO systems
proposed in \cite{Pan2016} is a variant of the block coordinate descent
(BCD) algorithm. Although this approach extends to MIMO systems, the
approximate problems solved in each iteration do not exhibit any convexity
and are thus not easy to solve, making the iterative algorithm not
suitable for practical implementation.

A low complexity algorithm is proposed in \cite{Zappone2017} to find
a KKT point of the nonconvex GEE optimization problem. The central
idea therein is to maximize in each iteration an approximate function
that is concave and a global lower bound of the original GEE function.
On the one hand, the maximum point of the concave approximate function
does not have a closed-form expression and it can only be found iteratively
by a general purpose optimization solver. On the other hand, this
sequential programming approach does not naturally extend to the SEE
problem because an approximate function that is a global lower bound
of the SEE function does not exist. This is also the case for the
GEE problem when the rate-dependent processing power consumption (due
to, e.g., coding and decoding, cf. \cite{Rost2010,Bjornson2015,Tervo2017})
is considered. Note that a global optimization technique is also proposed
in \cite{Zappone2017} for the GEE and SEE maximization problems,
which may serve as a benchmark in small problem instances only due
to the exponential complexity.

An iterative algorithm is proposed in \cite{He2014} to maximize the
SEE in MIMO systems (without QoS constraints). However, it has two
limitations. Firstly, it is a two layer algorithm for which the inner
layer consists of a BCD type algorithm which suffers from a high complexity
and a slow convergence rate. Secondly, only convergence in function
value is established and the convergence to a stationary point is
still left open. Besides, the algorithm is not applicable when the
rate-dependent processing power consumption is considered.

\begin{table}[t]
\center%
\begin{tabular}{|c|c|c|c|c|c|c|}
\hline
paper & D1 & D2 & D3 & GEE & SEE & QoS\tabularnewline
\hline
Tervo et al. \cite{Tervo2015} & \multirow{1}{*}{\checkmark} & \multirow{1}{*}{} & \multirow{1}{*}{} & \multirow{1}{*}{\checkmark} & \multirow{1}{*}{} & \checkmark\tabularnewline
\hline
Tervo et al. \cite{Tervo2017} & \checkmark & \checkmark &  & \checkmark & \checkmark & \checkmark\tabularnewline
\hline
Pan et al. \cite{Pan2016} & \checkmark &  & \checkmark &  & \checkmark & \checkmark\tabularnewline
\hline
Zappone et al. \cite{Zappone2017} & \checkmark &  & \checkmark & \checkmark &  & \checkmark\tabularnewline
\hline
He et al. \cite{He2014} & \checkmark &  & \checkmark &  & \checkmark & \tabularnewline
\hline
This paper & \checkmark & \checkmark & \checkmark & \checkmark & \checkmark & \checkmark\tabularnewline
\hline
\end{tabular}\caption{\label{tab:Summary-of-related-work}Summary of related work: addressed
difficulties}

\vspace{-1em}
\end{table}

\emph{Contributions. }In this paper, we study the GEE and SEE optimization
problems in multi-cell MIMO interference channels and propose novel
iterative algorithms that address the practical difficulties (D1)-(D3)
(see Table \ref{tab:Summary-of-related-work}), first without and
then with per-link QoS constraints; see Table \ref{tab:Summary-of-related-work}
for a comparison with some of the related works discussed above. The
proposed algorithms have the following attractive features:
\begin{itemize}
\item fast convergence as the structure of the original optimization problem
is preserved as much as possible in the approximate problem solved
in each iteration;
\item low complexity as each approximate problem is suitable for parallel
computation and its solution has a closed-form expression;
\item guaranteed convergence to a stationary point or a Karush-Kuhn-Tucker
(KKT) point.
\end{itemize}
The proposed algorithms are based on an extension of the recently
developed successive pseudoconvex approximation framework \cite{Yang_ConvexApprox}.
In each iteration, an approximate problem is solved, and the approximate
problem only needs to exhibit a weak form of convexity, namely, pseudoconvexity.
Among others, pseudoconvex optimization problems have two notable
properties: firstly, some special cases of pseudoconvex objective
functions (e.g., the ratio of positive convex and concave functions)
can be easily optimized and every stationary point is globally optimal,
and secondly, any direction pointing to an optimal point of a pseudoconvex
optimization problem is a descent direction of the objective function
(this property holds for convex functions as well but not for quasiconvex
functions). While the first property has been recognized and exploited
under the framework of fractional programming in many existing works
(see \cite{Ng2012,Ng2012a,Zappone2015,Zappone2017} and the references
therein), the second property has largely been overlooked. In this
paper, we argue that it plays a fundamental role in designing novel
iterative algorithms with provable convergence by showing repeatedly
that it paves the way to define an approximate problem that preserves
as much structure available in the original EE function as possible,
e.g., the partial concavity (convexity) in the numerator (denominator)
function and the division operator. Therefore, the proposed algorithm
presents a fast convergence behavior and enjoys an easy implementation.

We mention for the completeness of this paper that another popular
design metric is to maximize the minimum EE among all links. This
problem has been studied in \cite{Yang_MaxMinEE} where no rate-dependent
processing power consumption is considered. Our method proposed in
this paper cannot be applied to maximize the minimum EE, because the
minimum EE is a nondifferentiable function. To our best knowledge,
the minimum EE maximization problem with rate-dependent processing
power consumption is still an open problem.

\emph{Paper structure.} The rest of the paper is organized as follows.
In Sec. \ref{sec:Problem-Model} we introduce the system model and
problem formulation. The novel iterative algorithms are proposed in
Sections \ref{sec:GlobalEE}-\ref{sec:SEE-w/QoS} for the following
four problems: GEE maximization without QoS constraints, SEE maximization
without QoS constraints, GEE maximization with QoS constraints, and
SEE maximization with QoS constraints. Numerical results are reported
in Section \ref{sec:Simulations} and the paper is concluded in Sec.
\ref{sec:Concluding-Remarks}.

\emph{Notation: }We use $x$, $\mathbf{x}$ and $\mathbf{X}$ to denote
a scalar, vector and matrix, respectively. We use $\mathbf{X}^{H}$
and $\mathbf{X}^{*}$ to denote the Hermitian of $\mathbf{X}$ and
the complex conjugate of $\mathbf{X}$, respectively. The inner product
of two matrices $\mathbf{X}$ and $\mathbf{Y}$ is defined as $\mathbf{X}\bullet\mathbf{Y}\triangleq\Re(\textrm{tr}(\mathbf{X}^{H}\mathbf{Y}))$.
The operator $[\mathbf{X}]^{+}$ returns the projection of $\mathbf{X}$
onto the cone of positive semidefinite matrices. The gradient of $f(\mathbf{X})$
with respect to $\mathbf{X}^{*}$ and $\mathbf{X}_{k}^{*}$ is denoted
as $\nabla_{\mathbf{Q}^{*}}f(\mathbf{X})$ and $\nabla_{\mathbf{Q}_{k}^{*}}f(\mathbf{X})$,
respectively. $\nabla_{\mathbf{Q}^{*}}f(\mathbf{X})$ and $\nabla f(\mathbf{X})$
($\nabla_{\mathbf{Q}_{k}^{*}}f(\mathbf{X})$ and $\nabla_{k}f(\mathbf{X})$)
are used interchangeably when there is no ambiguity. When there are
multiple matrix variables $\mathbf{X}_{1},\mathbf{X}_{2},\ldots,\mathbf{X}_{K}$,
we use $\mathbf{X}$ as a compact notation to denote all of them:
$\mathbf{X}\triangleq(\mathbf{X}_{k})_{k=1}^{K}$. We also use $\mathbf{X}_{-k}$
to denote all matrix variables except $\mathbf{X}_{k}$: $\mathbf{X}_{-k}\triangleq(\mathbf{X}_{j})_{j=1,j\neq k}^{K}$.
The notation $\mathbf{0}\preceq\mathbf{X}\perp\mathbf{Y}\succeq\mathbf{0}$
denotes that $\mathbf{X}\succeq\mathbf{0}$, $\mathbf{Y}\succeq\mathbf{0}$
and $\Re(\textrm{tr}(\mathbf{X}^{H}\mathbf{Y}))=0$. Similarly $\mathbf{0}\leq\mathbf{x}\perp\mathbf{y}\geq\mathbf{0}$
denotes that $\mathbf{x}\geq\mathbf{0}$, $\mathbf{y}\geq\mathbf{0}$
and $\mathbf{x}^{H}\mathbf{y}=0$.

\section{\label{sec:Problem-Model}System Model and Problem Formulation}

\begin{figure}[t]
\center\includegraphics[bb=80bp 320bp 420bp 670bp,clip,scale=0.6]{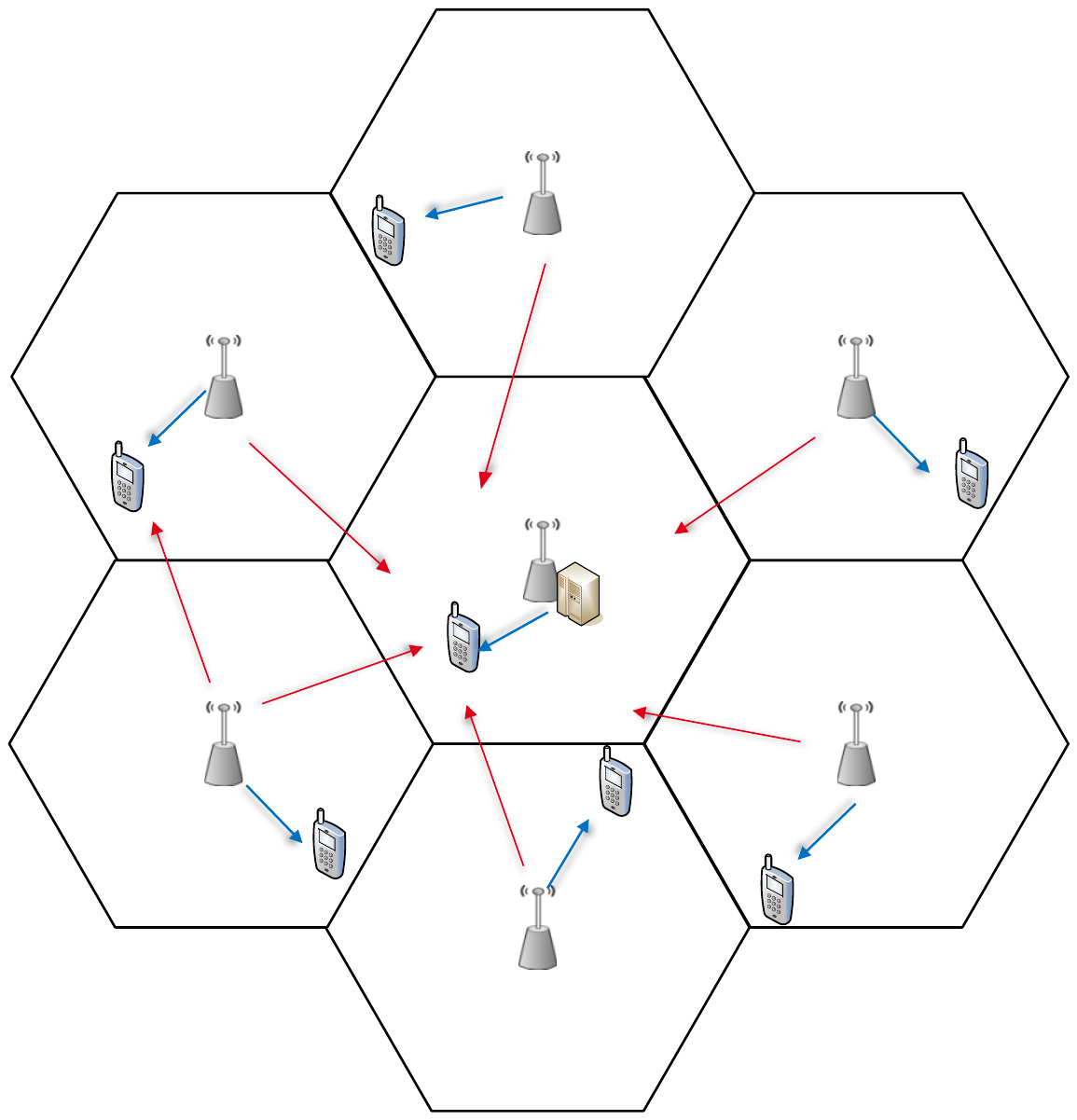}

\caption{\label{fig:System-topology}System topology of 1-tiered small cell
interferers (with a central processing unit in CRAN). The BSs and
users may be equipped with multiple antennas.}

\vspace{-1em}
\end{figure}

We consider a downlink MIMO multi-cell system as depicted in Figure
\ref{fig:System-topology}, where the number of cells is $K$. We
assume for simplicity that each cell is serving one user, but the
results can be generalized to the case that each cell is serving multiple
users. The number of transmit antennas at the BS of cell $k$ is $M_{k}$,
and the number of receive antennas of user $k$ served by cell $k$
is $N_{k}$. We denote $\mathbf{H}_{kk}$ as the channel matrix from
BS $k$ to user $k$, and $\mathbf{H}_{kj}$ as the downlink channel
matrix from BS $j$ to user $k$. We assume that all $K$ users are
active and the multi-user interference is treated as noise, so the
downlink transmission rate of the $k$-th user is:
\begin{align}
r_{k}(\mathbf{Q}_{k},\mathbf{Q}_{-k})\triangleq & \log\det\left(\mathbf{I}+\mathbf{R}_{k}(\mathbf{Q}_{-k})^{-1}\mathbf{H}_{kk}\mathbf{Q}_{k}\mathbf{H}_{kk}^{H}\right),\label{eq:rate}
\end{align}
where $\mathbf{Q}_{k}\triangleq\mathbb{E}\left[\mathbf{x}_{k}\mathbf{x}_{k}^{H}\right]$
is BS $k$'s transmit covariance matrix, $\mathbf{Q}_{-k}$ is a compact
notation denoting all transmit covariance matrices except $\mathbf{Q}_{k}$:
$\mathbf{Q}_{-k}=(\mathbf{Q}_{j})_{j\neq k}$, and $\mathbf{R}_{k}(\mathbf{Q}_{-k})\triangleq\sigma_{k}^{2}\mathbf{I}+\sum_{j\neq k}\mathbf{H}_{kj}\mathbf{Q}_{j}\mathbf{H}_{kj}^{H}$
is the noise plus interference covariance matrix experienced by user
$k$.

The power consumption at BS $k$ can be approximated by the following
equation:
\[
p_{k}(\mathbf{Q})\triangleq P_{0,k}+\rho_{k}\textrm{tr}(\mathbf{Q}_{k})+g_{k}(r_{k}(\mathbf{Q})),
\]
where $\mathbf{Q}\triangleq(\mathbf{Q}_{k})_{k=1}^{K}$, $P_{0,k}$
is the circuit power consumption and $\rho_{k}\geq1$ is the inverse
of the power amplifier efficiency at the transmitter, while $g_{k}(x)$
is a monotonic and differentiable function of $x$ with $g_{k}(0)=0$,
which reflects the rate-dependent processing power consumption, e.g.,
required for coding and decoding \cite{Rost2010,Bjornson2015,Tervo2017}.
The typical values of $P_{0,k}$ and $\rho_{k}$ depend on the types
of the cell, e.g., macro cell, remote radio head, and micro cell.
Interested readers are referred to \cite[Table 8]{EarthD2.3} for
its typical values.

Optimizing the EE of multiple links in the network simultaneously
is a typical multi-objective optimization problem, which could be
modeled in several different ways. For example, the GEE, which is
defined as the ratio between the sum transmission rate and the total
consumed power, models the EE of the whole network:
\begin{align}
\underset{\mathbf{Q}}{\textrm{maximize}}\quad & f_{G}(\mathbf{Q})\triangleq\frac{\sum_{k=1}^{K}r_{k}(\mathbf{Q})}{\sum_{k=1}^{K}p_{k}(\mathbf{Q})}\nonumber \\
\textrm{subject to}\quad & \mathbf{Q}_{k}\succeq\mathbf{0},\,\textrm{tr}(\mathbf{Q}_{k})\leq P_{k},\,\forall k,\label{eq:original problem-GlobalEE}
\end{align}
where $P_{k}$ is BS $k$'s (predefined) sum transmission power budget
and the optimization variable is the (complex-valued) transmit covariance
matrices $\mathbf{Q}=(\mathbf{Q}_{k})_{k=1}^{K}$.

To maximize the GEE, users suffering from bad channel conditions may
not be able to transmit, because increasing the transmit power in
the denominator may not lead to a notable increase in the transmission
rate in the numerator. Another popular design approach is the so-called
socially optimal approach, which aims at finding the optimal point
that maximizes the sum EE (SEE) over all users:
\begin{align}
\underset{\mathbf{Q}}{\textrm{maximize}}\quad & f_{S}(\mathbf{Q})\triangleq\sum_{k=1}^{K}\frac{r_{k}(\mathbf{Q})}{p_{k}(\mathbf{Q})}\nonumber \\
\textrm{subject to}\quad & \mathbf{Q}_{k}\succeq\mathbf{0},\,\textrm{tr}(\mathbf{Q}_{k})\leq P_{k},\,\forall k.\label{eq:original-problem-SumEE}
\end{align}
Note that the objective function $f_{S}(\mathbf{Q})$ in (\ref{eq:original-problem-SumEE})
is a sum of multiple fractional functions, each is the ratio of the
nonconcave function $r_{k}(\mathbf{Q})$ and the nonconvex function
$p_{k}(\mathbf{Q})$.

In the previous formulations (\ref{eq:original problem-GlobalEE})
and (\ref{eq:original-problem-SumEE}), there are no QoS constraints
specifying each link's minimum guaranteed transmission rate. To incorporate
the QoS constraints into the EE optimization problems, we modify the
GEE optimization problem (\ref{eq:original problem-GlobalEE}) as
follows:
\begin{align}
\underset{\mathbf{Q}}{\textrm{maximize}}\quad & f_{G}(\mathbf{Q})=\frac{\sum_{k=1}^{K}r_{k}(\mathbf{Q})}{\sum_{k=1}^{K}p_{k}(\mathbf{Q})}\nonumber \\
\textrm{subject to}\quad & \mathbf{Q}_{k}\succeq\mathbf{0},\,\textrm{tr}(\mathbf{Q}_{k})\leq P_{k},r_{k}(\mathbf{Q})\geq R_{k},\,\forall k,\label{eq:original-problem-GlobalEE-w/QoS}
\end{align}
and we assume that the solution set of (\ref{eq:original-problem-GlobalEE-w/QoS})
is nonempty. In contrast to problem (\ref{eq:original problem-GlobalEE}),
problem (\ref{eq:original-problem-GlobalEE-w/QoS}) has a nonconvex
constraint set due to the nonconvex QoS constraints and is thus more
challenging. Therefore we study (\ref{eq:original problem-GlobalEE})
and (\ref{eq:original-problem-GlobalEE-w/QoS}) separately.

Similarly, the SEE optimization problem subject to QoS constraints
is modeled as follows:
\begin{align}
\underset{\mathbf{Q}}{\textrm{maximize}}\quad & f_{S}(\mathbf{Q})=\sum_{k=1}^{K}\frac{r_{k}(\mathbf{Q})}{p_{k}(\mathbf{Q})}\nonumber \\
\textrm{subject to}\quad & \mathbf{Q}_{k}\succeq\mathbf{0},\,\textrm{tr}(\mathbf{Q}_{k})\leq P_{k},r_{k}(\mathbf{Q})\geq R_{k},\,\forall k.\label{eq:original-problem-SumEE-w/QoS}
\end{align}
In Sections \ref{sec:GlobalEE}-\ref{sec:SEE-w/QoS}, we propose novel
iterative algorithms that can efficiently find a stationary/KKT point
of problems (\ref{eq:original problem-GlobalEE})-(\ref{eq:original-problem-SumEE-w/QoS}).

\section{\label{sec:GlobalEE}The Proposed Iterative Algorithm For Global
Energy Efficiency Maximization}

To design an iterative algorithm for problem (\ref{eq:original problem-GlobalEE})
that enjoys a low complexity but at the same time a fast convergence
behavior, we need on the one hand to address the issue of the nonconvexity
in the objective function, and, on the other hand, to preserve the
original problem's structure as much as possible. Towards this end,
we propose an iterative algorithm based on the successive pseudoconvex
approximation framework developed in \cite{Yang_ConvexApprox}.

To start with, we introduce the definition of pseudoconvex functions:
a function $f(\mathbf{x})$ is said to be \emph{pseudoconvex} if \cite{Mangasarian_NonlinearProgramming}
\begin{equation}
f(\mathbf{y})<f(\mathbf{x})\Longrightarrow(\mathbf{y}-\mathbf{x})^{T}\nabla f(\mathbf{x})<0.\label{eq:def-pseudo-convex-function}
\end{equation}
In other words, $f(\mathbf{y})<f(\mathbf{x})$ implies $\mathbf{y}-\mathbf{x}$
is a \emph{descent direction} of $f(\mathbf{x})$ \cite{bertsekas1999nonlinear}.
A function $f(\mathbf{x})$ is pseudoconcave if $-f(\mathbf{x})$
is pseudoconvex. We remark that the (strong) convexity of a function
implies that the function is pseudoconvex, which in turn implies that
the function is quasiconvex, but the reverse is generally not true;
see \cite[Figure 1]{Yang_ConvexApprox}.

The proposed iterative algorithm for problem (\ref{eq:original problem-GlobalEE})
consists of solving a sequence of successively refined approximate
problems. In iteration $t$, the approximate problem defined around
a given point $\mathbf{Q}^{t}$ consists of maximizing an approximate
function, denoted as $\tilde{f}(\mathbf{Q};\mathbf{Q}^{t})$, under
the same constraints as (\ref{eq:original problem-GlobalEE}). The
lack of concavity in the objective function should be properly compensated
so that the approximate problems are much easier to solve than the
original problem (\ref{eq:original problem-GlobalEE}).

The numerator functions $(r_{k}(\mathbf{Q}))_{k=1}^{K}$ are not concave
and the denominator functions $(p_{k}(\mathbf{Q}))_{k=1}^{K}$ are
not convex in $\mathbf{Q}$. Meanwhile, the function $r_{k}(\mathbf{Q})$
is concave in component $\mathbf{Q}_{k}$, and the function $P_{0,k}+\rho_{k}\textrm{tr}(\mathbf{Q}_{k})$
in $p_{k}(\mathbf{Q})$ is convex in component $\mathbf{Q}_{k}$.
Exploiting this partial concavity may notably accelerate the convergence,
as shown in \cite{Scutarib} and other works. Therefore, we approximate
the numerator function $\sum_{j=1}^{K}r_{j}(\mathbf{Q})$ with respect
to (w.r.t.) $\mathbf{Q}_{k}$ at the point $\mathbf{Q}^{t}$ by a
function denoted as $\tilde{r}_{G,k}(\mathbf{Q}_{k};\mathbf{Q}^{t})$,
which is obtained by fixing the other variables $\mathbf{Q}_{-k}$
in $r_{k}(\mathbf{Q}_{k},\mathbf{Q}_{-k})$ and linearizing only the
functions $\{r_{j}(\mathbf{Q})\}_{j\neq k}$ that are not concave
in $\mathbf{Q}_{k}$:
\begin{equation}
\tilde{r}_{G,k}(\mathbf{Q}_{k};\mathbf{Q}^{t})\triangleq r_{k}(\mathbf{Q}_{k},\mathbf{Q}_{-k}^{t})+{\textstyle \sum_{j\neq k}}(\mathbf{Q}_{k}-\mathbf{Q}_{k}^{t})\bullet\nabla_{k}r_{j}(\mathbf{Q}^{t}),\label{eq:approximate-r}
\end{equation}
where $\mathbf{X}\bullet\mathbf{Y}\triangleq\Re(\textrm{tr}(\mathbf{X}^{H}\mathbf{Y}))$
and $\nabla_{k}r_{j}(\mathbf{Q})$ is the Jacobian matrix of $r_{j}(\mathbf{Q})$
with respect to $\mathbf{Q}_{k}^{*}$ (the complex conjugate of $\mathbf{Q}_{k}$).
Since $\tilde{r}_{G,k}(\mathbf{Q}_{k};\mathbf{Q}^{t})$ is concave
in $\mathbf{Q}_{k}$, $\sum_{k=1}^{K}\tilde{r}_{G,k}(\mathbf{Q}_{k};\mathbf{Q}^{t})$
is concave in $\mathbf{Q}$. Similarly, we approximate the denominator
function $p_{k}(\mathbf{Q})$ by a convex function $\tilde{p}_{G,k}(\mathbf{Q}_{k};\mathbf{Q}^{t})$
which is obtained by keeping $P_{0,k}+\rho_{k}\textrm{tr}(\mathbf{Q}_{k})$
and linearizing the nonconvex part $p_{k}(\mathbf{Q})$ w.r.t. $\mathbf{Q}_{k}$
at the point $\mathbf{Q}=\mathbf{Q}^{t}$:
\begin{align}
\tilde{p}_{G,k}(\mathbf{Q}_{k};\mathbf{Q}^{t})\triangleq\; & P_{0,k}+\rho_{k}\textrm{tr}(\mathbf{Q}_{k})+g_{k}(r_{k}(\mathbf{Q}^{t}))\nonumber \\
 & +{\textstyle \sum_{j=1}^{K}}(\mathbf{Q}_{k}-\mathbf{Q}_{k}^{t})\bullet\nabla_{k}g_{j}(r_{j}(\mathbf{Q}^{t})),\label{eq:approximate-p}
\end{align}
and $\tilde{p}_{G,k}(\mathbf{Q}_{k};\mathbf{Q}^{t})$ is positive
and convex. This paves the way to define the following approximate
function of the original objective function $f(\mathbf{Q})$ at point
$\mathbf{Q}^{t}$, denoted as $\tilde{f}_{G}(\mathbf{Q};\mathbf{Q}^{t})$:
\begin{equation}
\tilde{f}_{G}(\mathbf{Q};\mathbf{Q}^{t})\triangleq\frac{\sum_{k=1}^{K}\tilde{r}_{G,k}(\mathbf{Q}_{k};\mathbf{Q}^{t})}{\sum_{k=1}^{K}\tilde{p}_{G,k}(\mathbf{Q}_{k};\mathbf{Q}^{t})},\label{eq:notation-approximate_function}
\end{equation}
The approximate function $\tilde{f}_{G}(\mathbf{Q};\mathbf{Q}^{t})$
has some important properties as we outline.

Firstly, the approximate function $\tilde{f}_{G}(\mathbf{Q};\mathbf{Q}^{t})$
is still nonconcave, but it is a fractional function of a nonnegative
concave function $\sum_{k=1}^{K}\tilde{r}_{G,k}(\mathbf{Q}_{k};\mathbf{Q}^{t})$
and a positive linear function $\sum_{k=1}^{K}\tilde{p}_{G,k}(\mathbf{Q}_{k};\mathbf{Q}^{t})$,
which is thus pseudoconcave \cite{Yang_ConvexApprox}.

Secondly, the approximate function $\tilde{f}_{G}(\mathbf{Q};\mathbf{Q}^{t})$
is differentiable and its gradient is the same as that of the original
function $f_{G}(\mathbf{Q})$ at the point $\mathbf{Q}^{t}$ where
the approximate function $\tilde{f}_{G}(\mathbf{Q};\mathbf{Q}^{t})$
is defined. To see this, we remark that $\left.\nabla_{\mathbf{Q}_{k}^{*}}\tilde{r}_{j}(\mathbf{Q}_{j};\mathbf{Q}^{t})\right|_{\mathbf{Q}=\mathbf{Q}^{t}}=\mathbf{0}$
if $j\neq k$, and\begin{subequations}\label{eq:equal-r}
\begin{align}
\left.\nabla_{k}\tilde{r}_{G,k}(\mathbf{Q}_{k};\mathbf{Q}^{t})\right|_{\mathbf{Q}=\mathbf{Q}^{t}} & =\left.\nabla_{k}\left({\textstyle \sum_{j=1}^{K}}r_{j}(\mathbf{Q})\right)\right|_{\mathbf{Q}=\mathbf{Q}^{t}},\label{eq:equal-r-gradient}\\
\tilde{r}_{G,k}(\mathbf{Q}_{k}^{t};\mathbf{Q}^{t}) & =r_{k}(\mathbf{Q}^{t}).\label{eq:equal-r-value}
\end{align}
\end{subequations}Similarly, $\nabla_{k}\tilde{g}_{G,j}(\mathbf{Q}_{j};\mathbf{Q}^{t})=\mathbf{0}$
if $j\neq k$ and\begin{subequations}\label{eq:equal-g}
\begin{align}
\left.\nabla_{k}\tilde{p}_{G,k}(\mathbf{Q}_{k};\mathbf{Q}^{t})\right|_{\mathbf{Q}=\mathbf{Q}^{t}} & =\left.\nabla_{k}\left({\textstyle \sum_{j=1}^{K}}p_{j}(\mathbf{Q})\right)\right|_{\mathbf{Q}=\mathbf{Q}^{t}},\label{eq:equal-g-gradient}\\
\tilde{p}_{G,k}(\mathbf{Q}_{k}^{t};\mathbf{Q}^{t}) & =p_{k}(\mathbf{Q}^{t}).\label{eq:equal-g-value}
\end{align}
\end{subequations} Based on the observations in (\ref{eq:equal-r})-(\ref{eq:equal-g}),
it can be verified that the gradient of the approximate function $\tilde{f}_{G}(\mathbf{Q};\mathbf{Q}^{t})$
is the same as that of the original function $f_{G}(\mathbf{Q})$
at the point $\mathbf{Q}^{t}$:
\begin{align}
 & \left.\nabla_{k}\tilde{f}_{G}(\mathbf{Q};\mathbf{Q}^{t})\right|_{\mathbf{Q}=\mathbf{Q}^{t}}\nonumber \\
=\, & \frac{\nabla_{k}\tilde{r}_{k}(\mathbf{Q}_{k}^{t};\mathbf{Q}^{t})}{\sum_{j=1}^{K}\tilde{p}_{j}(\mathbf{Q}_{j}^{t};\mathbf{Q}^{t})}-\frac{(\sum_{j=1}^{K}\tilde{r}_{j}(\mathbf{Q}_{j}^{t};\mathbf{Q}^{t}))\nabla_{k}\tilde{p}_{k}(\mathbf{Q}_{k}^{t};\mathbf{Q}^{t})}{(\sum_{j=1}^{K}\tilde{p}_{j}(\mathbf{Q}_{j}^{t};\mathbf{Q}^{t}))^{2}}\nonumber \\
=\, & \frac{\nabla_{k}(\sum_{j=1}^{K}r_{j}(\mathbf{Q}^{t}))}{\sum_{j=1}^{K}p_{j}(\mathbf{Q}^{t})}-\frac{(\sum_{j=1}^{K}r_{j}(\mathbf{Q}^{t}))\nabla_{k}(\sum_{j=1}^{K}p_{j}(\mathbf{Q}^{t}))}{(\sum_{j=1}^{K}p_{k}(\mathbf{Q}^{t}))^{2}}\nonumber \\
=\, & \left.\nabla_{k}f_{G}(\mathbf{Q})\right|_{\mathbf{Q}=\mathbf{Q}^{t}},\label{eq:euqal-gradient-f}
\end{align}
where the first and third equality is the expression of $\nabla_{k}\tilde{f}_{G}(\mathbf{Q};\mathbf{Q}^{t})$
and $\nabla_{k}f_{G}(\mathbf{Q})$, respectively, and the second equality
follows from (\ref{eq:equal-r})-(\ref{eq:equal-g}).

At iteration $t$ of the proposed algorithm, the approximate problem
defined at the point $\mathbf{Q}^{t}$ is to maximize the approximate
function $\tilde{f}_{G}(\mathbf{Q};\mathbf{Q}^{t})$ defined in (\ref{eq:notation-approximate_function})
subject to the same constraints as in the original problem (\ref{eq:original problem-GlobalEE}):\begin{subequations}\label{eq:approximate-problem}
\begin{align}
\underset{\mathbf{Q}}{\textrm{maximize}}\quad & \tilde{f}_{G}(\mathbf{Q};\mathbf{Q}^{t})\nonumber \\
\textrm{subject to}\quad & \mathbf{Q}_{k}\succeq\mathbf{0},\textrm{tr}(\mathbf{Q}_{k})\leq P_{k},\,k=1,\ldots,K.\label{eq:approximate-problem-formulation}
\end{align}
and its (globally) optimal point is denoted as $\mathbb{B}\mathbf{Q}^{t}$:
\begin{equation}
\mathbb{B}\mathbf{Q}^{t}\triangleq\underset{(\mathbf{Q}_{k}\succeq\mathbf{0},\textrm{tr}(\mathbf{Q}_{k})\leq P_{k})_{k=1}^{K}}{\arg\max}\;\tilde{f}_{G}(\mathbf{Q};\mathbf{Q}^{t}).\label{eq:approximate-problem-optimal-point}
\end{equation}
\end{subequations}Since problem (\ref{eq:approximate-problem-formulation})
is pseudoconvex, all of its stationary points are globally optimal
\cite[Th. 9.3.3]{Mangasarian_NonlinearProgramming}. As we will show
shortly, $\mathbb{B}\mathbf{Q}^{t}$ is unique.

Due to the above mentioned pseudoconcavity, differentiability and
equal gradient condition (\ref{eq:euqal-gradient-f}) at $\mathbf{Q}^{t}$
of the approximate function $\tilde{f}(\mathbf{Q};\mathbf{Q}^{t})$
defined in (\ref{eq:notation-approximate_function}), solving the
approximate problem (\ref{eq:approximate-problem}) yields an ascent
direction of the original objective function $f_{G}(\mathbf{Q})$
at $\mathbf{Q}^{t}$, unless $\mathbf{Q}^{t}$ is already a stationary
point of problem (\ref{eq:original problem-GlobalEE}), as stated
in the following proposition.
\begin{prop}
[Stationary point and ascent direction]\label{prop:ascent-direction}A
point $\mathbf{Q}^{t}$ is a stationary point of (\ref{eq:original problem-GlobalEE})
if and only if $\mathbf{Q}^{t}=\mathbb{B}\mathbf{Q}^{t}$. If $\mathbf{Q}^{t}$
is not a stationary point of (\ref{eq:original problem-GlobalEE}),
then $\mathbb{B}\mathbf{Q}^{t}-\mathbf{Q}^{t}$ is an ascent direction
of $f_{G}(\mathbf{Q})$ in the sense that
\[
(\mathbb{B}\mathbf{Q}^{t}-\mathbf{Q}^{t})\bullet\nabla f_{G}(\mathbf{Q}^{t})>0.
\]
\end{prop}
\begin{IEEEproof}
By (\ref{eq:approximate-problem-optimal-point}), if $\mathbf{Q}^{t}=\mathbb{B}\mathbf{Q}^{t}$,
then $\mathbf{Q}^{t}$ is an optimal point of the following problem:
\[
\mathbf{Q}^{t}=\underset{(\mathbf{Q}_{k}\succeq\mathbf{0},\textrm{tr}(\mathbf{Q}_{k})\leq P_{k})_{k=1}^{K}}{\arg\max}\;\tilde{f}_{G}(\mathbf{Q};\mathbf{Q}^{t}).
\]
According to the first-order optimality condition, the following inequality
is satisfied:
\[
(\mathbf{Q}-\mathbf{Q}^{t})\bullet\nabla\tilde{f}_{G}(\mathbf{Q}^{t};\mathbf{Q}^{t})\geq0
\]
for all $\mathbf{Q}$ such that $\mathbf{Q}_{k}\succeq\mathbf{0}$
and $\textrm{tr}(\mathbf{Q}_{k})\leq P_{k}$, $k=1,\ldots,K$. Since
$\nabla\tilde{f}_{G}(\mathbf{Q}^{t};\mathbf{Q}^{t})=\nabla f_{G}(\mathbf{Q}^{t})$,
cf. (\ref{eq:euqal-gradient-f}), the above inequality is
\[
(\mathbf{Q}-\mathbf{Q}^{t})\bullet\nabla f_{G}(\mathbf{Q}^{t})\geq0
\]
for all $\mathbf{Q}$ such that $\mathbf{Q}_{k}\succeq\mathbf{0}$
and $\textrm{tr}(\mathbf{Q}_{k})\leq P_{k}$, $k=1,\ldots,K$. This
is the first order optimality condition of the original problem (\ref{eq:original problem-GlobalEE})
and $\mathbf{Q}^{t}$ is thus a stationary point of (\ref{eq:original problem-GlobalEE}).

If $\mathbf{Q}^{t}\neq\mathbb{B}\mathbf{Q}^{t}$, then
\[
\tilde{f}_{G}(\mathbb{B}\mathbf{Q}^{t};\mathbf{Q}^{t})>\tilde{f}_{G}(\mathbf{Q}^{t};\mathbf{Q}^{t}).
\]
Since $\tilde{f}_{G}(\mathbf{Q};\mathbf{Q}^{t})$ is a pseudoconcave
function in $\mathbf{Q}$, it follows from the definition (\ref{eq:def-pseudo-convex-function})
that
\[
0<(\mathbb{B}\mathbf{Q}^{t}-\mathbf{Q}^{t})\bullet\nabla\tilde{f}_{G}(\mathbf{Q}^{t};\mathbf{Q}^{t})=(\mathbb{B}\mathbf{Q}^{t}-\mathbf{Q}^{t})\bullet\nabla f_{G}(\mathbf{Q}^{t}),
\]
where the equality comes from the fact that $\nabla\tilde{f}_{G}(\mathbf{Q}^{t};\mathbf{Q}^{t})=\nabla f_{G}(\mathbf{Q}^{t})$.
The proof is thus completed.
\end{IEEEproof}
Since $\mathbb{B}\mathbf{Q}^{t}-\mathbf{Q}^{t}$ is an ascent direction
of $f_{G}(\mathbf{Q})$ at $\mathbf{Q}=\mathbf{Q}^{t}$ according
to Proposition \ref{prop:ascent-direction}, there exists a scalar
$\gamma^{t}\in(0,1]$ such that $f_{G}(\mathbf{Q}^{t}+\gamma^{t}(\mathbb{B}\mathbf{Q}^{t}-\mathbf{Q}^{t}))>f_{G}(\mathbf{Q}^{t})$
\cite[8.2.1]{Ortega&Rheinboldt}. In practice, the stepsize $\gamma^{t}$
is usually obtained by either the exact line search or the successive
line search. Performing the exact line search consists of solving
an optimization problem
\[
\max_{0\leq\gamma\leq1}f(\mathbf{Q}^{t}+\gamma(\mathbb{B}\mathbf{Q}^{t}-\mathbf{Q}^{t})).
\]
Since the objective function $f(\mathbf{Q})$ is nonconcave, the above
optimization problem is nonconvex and not trivial to solve. Therefore,
we adopt the successive line search to calculate the stepsize $\gamma^{t}$.
That is, given two scalars $0<\alpha<1$ and $0<\beta<1$, $\gamma^{t}$
is set to be $\gamma^{t}=\beta^{m_{t}}$, where $m_{t}$ is the smallest
nonnegative integer $m$ satisfying the following inequality:
\begin{align}
f_{G}(\mathbf{Q}^{t}+\beta^{m}(\mathbb{B}\mathbf{Q}^{t}-\mathbf{Q}^{t})) & \geq\nonumber \\
f_{G}(\mathbf{Q}^{t})+\alpha\beta^{m} & \nabla f_{G}(\mathbf{Q}^{t})\bullet(\mathbb{B}\mathbf{Q}^{t}-\mathbf{Q}^{t}).\label{eq:line-search}
\end{align}
Note that the successive line search is carried out over the original
objective function $f(\mathbf{Q})$ defined in (\ref{eq:original problem-GlobalEE}).

After the stepsize $\gamma^{t}$ is found, the variable $\mathbf{Q}$
is updated as
\begin{equation}
\mathbf{Q}^{t+1}=\mathbf{Q}^{t}+\gamma^{t}(\mathbb{B}\mathbf{Q}^{t}-\mathbf{Q}^{t}).\label{eq:variable-update}
\end{equation}
The resulting sequence $\{f_{G}(\mathbf{Q}^{t})\}_{t}$ is increasing:
\begin{align*}
f_{G}(\mathbf{Q}^{t+1}) & =f_{G}(\mathbf{Q}^{t}+\beta^{m_{t}}(\mathbb{B}\mathbf{Q}^{t}-\mathbf{Q}^{t}))\\
 & \geq f_{G}(\mathbf{Q}^{t})+\alpha\beta^{m_{t}}\nabla f_{G}(\mathbf{Q}^{t})\bullet(\mathbb{B}\mathbf{Q}^{t}-\mathbf{Q}^{t})\\
 & \geq f_{G}(\mathbf{Q}^{t}),\;\forall t,
\end{align*}
where the first and second inequality comes from the definition of
the successive line search (\ref{eq:line-search}) and Proposition
\ref{prop:ascent-direction}, respectively.

The proposed algorithm is formally summarized in Algorithm \ref{alg-GlobalEE}
and its convergence properties are given in the following theorem.

\begin{algorithm}[t]
\textbf{S0:} $\mathbf{Q}^{0}=\mathbf{0}$, $t=0$, and a stopping
criterion $\varepsilon$.

\textbf{S1: }Compute $\mathbb{B}\mathbf{Q}^{t}$ by solving problem
(\ref{eq:approximate-problem}):

$\qquad$\textbf{S1.0: }$s^{t,0}=0$, $\tau=0$, and a stopping criterion
$\epsilon$.

$\qquad$\textbf{S1.1: }Compute $\mathbf{Q}_{k}^{\star}(s^{t,\tau})$
by (\ref{eq:dinkelbach-parallel-subproblem}).

$\qquad$\textbf{S1.2: }Compute $s^{t,\tau+1}$ by (\ref{eq:dinkelbach-alpha-parallel}).

$\qquad$\textbf{S1.3: }If $|s^{t,\tau+1}-s^{t,\tau}|<\epsilon$,
then $\mathbb{B}\mathbf{Q}^{t}=\mathbf{Q}^{\star}(s^{t,\tau})$. $\qquad$$\qquad$$\;$$\;$$\,$Otherwise
$\tau\leftarrow\tau+1$ and go to \textbf{S1.1}.

\textbf{S2: }Compute $\gamma^{t}$ by the successive line search (\ref{eq:line-search}).

\textbf{S3: }Update $\mathbf{Q}^{t+1}$ according to (\ref{eq:variable-update}).

\textbf{S4: }If $\left\Vert \mathbb{B}\mathbf{Q}^{t}-\mathbf{Q}^{t}\right\Vert \leq\varepsilon$,
then STOP; otherwise $t\leftarrow t+1$ and $\quad$$\;$$\,$$\,$go
to \textbf{S1}.

\caption{\label{alg-GlobalEE}The successive pseudoconvex approximation method
for GEE maximization (\ref{eq:original problem-GlobalEE})}
\end{algorithm}
\begin{thm}
[Convergence to a stationary point]\label{thm:convergence}The sequence
$\{\mathbf{Q}^{t}\}$ generated by Algorithm \ref{alg-GlobalEE} has
a limit point, and every limit point is a stationary point of problem
(\ref{eq:original problem-GlobalEE}).
\end{thm}
\begin{IEEEproof}
The constraint set of problem (\ref{eq:original problem-GlobalEE}),
namely, $\{(\mathbf{Q}_{k})_{k=1}^{K}:\mathbf{Q}_{k}\succeq\mathbf{0},\textrm{tr}(\mathbf{Q}_{k})\leq P_{k}\}$,
is nonempty and bounded. The sequence $\{\mathbf{Q}^{t}\}_{t}$ is
thus bounded and has a limit point. Then the latter statement can
be proved following the same line of analysis as \cite[Theorem 1]{Yang_ConvexApprox}
and is thus not duplicated here.
\end{IEEEproof}
In Step 1 of Algorithm \ref{alg-GlobalEE}, a constrained pseudoconvex
optimization problem, namely, problem (\ref{eq:approximate-problem})
must be solved. Since the optimal point $\mathbb{B}\mathbf{Q}^{t}$
does not have a closed-form expression, we apply the Dinkelbach's
algorithm \cite{Zappone2015} to solve problem (\ref{eq:approximate-problem})
iteratively: at iteration $\tau$ of Dinkelbach's algorithm, the following
problem is solved for a given and fixed $s^{t,\tau}$ ($s^{t,0}$
can be set to 0):
\begin{align}
\underset{\mathbf{Q}}{\textrm{maximize}}\quad & {\textstyle \sum_{k=1}^{K}}\tilde{r}_{G,k}(\mathbf{Q}_{k};\mathbf{Q}^{t})-s^{t,\tau}{\textstyle \sum_{k=1}^{K}\tilde{p}_{G,k}(\mathbf{Q}_{k};\mathbf{Q}^{t})}\nonumber \\
\textrm{subject to}\quad & \mathbf{Q}_{k}\succeq\mathbf{0},\,\textrm{tr}(\mathbf{Q}_{k})\leq P_{k},\;\forall k.\label{eq:dinkelbach-parallel}
\end{align}
Since problem (\ref{eq:dinkelbach-parallel}) is well decoupled across
different variables, it can be decomposed component-wise into many
smaller optimization problems that can be solved in parallel: for
all $k=1,\ldots,K$,\begin{subequations}\label{eq:dinkelbach-parallel-subproblem}
\begin{align}
\underset{\mathbf{Q}}{\textrm{maximize}}\quad & \tilde{r}_{G,k}(\mathbf{Q}_{k};\mathbf{Q}^{t})-s^{t,\tau}\tilde{p}_{G,k}(\mathbf{Q}_{k};\mathbf{Q}^{t})\nonumber \\
\textrm{subject to}\quad & \mathbf{Q}_{k}\succeq\mathbf{0},\,\textrm{tr}(\mathbf{Q}_{k})\leq P_{k}.\label{eq:dinkelbach-parallel-subproblem-formulation}
\end{align}
This problem is convex and its (unique) optimal point has a closed-form
expression based on the generalized waterfilling solution \cite[Lemma 2]{Yang_stochastic}:
\begin{align}
\mathbf{Q}_{k}^{\star}(s^{t,\tau}) & \triangleq\underset{\mathbf{Q}_{k}\succeq\mathbf{0},\textrm{tr}(\mathbf{Q}_{k})\leq P_{k}}{\arg\max}\left\{ \tilde{r}_{G,k}(\mathbf{Q}_{k};\mathbf{Q}^{t})-s^{t,\tau}\tilde{p}_{G,k}(\mathbf{Q}_{k};\mathbf{Q}^{t})\right\} \nonumber \\
 & =\mathbf{V}[\mathbf{I}-\boldsymbol{\Sigma}^{-1}]^{+}\mathbf{V}^{H},\label{eq:dinkelbach-parallel-subproblem-waterfilling}
\end{align}
\end{subequations}where $[\mathbf{X}]^{+}$ denotes the projection
of $\mathbf{X}$ onto the cone of positive semidefinite matrices,
$(\mathbf{V},\boldsymbol{\Sigma})$ is the generalized eigenvalue
decomposition of $(\mathbf{H}_{kk}^{H}\mathbf{R}_{k}(\mathbf{Q}_{-k}^{t})^{-1}\mathbf{H}_{kk},(s^{t,\tau}\rho_{k}+\mu^{\star})\mathbf{I}+s^{t,\tau}\nabla_{\mathbf{Q}_{k}^{*}}(\sum_{j=1}^{K}g_{j}(r_{j}(\mathbf{Q}^{t})))-\sum_{j\neq k}\nabla_{k}r_{j}(\mathbf{Q}^{t}))$,
and $\mu^{\star}$ is the Lagrange multiplier such that $0\leq\mu^{\star}\perp\textrm{tr}(\mathbf{Q}_{k}^{\star}(s^{t,\tau}))-P_{k}\leq0$,
which can easily be found by bisection.

After $(\mathbf{Q}_{k}^{\star}(s^{t,\tau}))_{k=1}^{K}$ is obtained,
$s^{t,\tau}$ is updated as follows:
\begin{equation}
s^{t,\tau+1}=\frac{\sum_{k=1}^{K}\tilde{r}_{G,k}(\mathbf{Q}_{k}^{\star}(s^{t,\tau});\mathbf{Q}^{t})}{\sum_{k=1}^{K}\tilde{p}_{G,k}(\mathbf{Q}^{\star}(s^{t,\tau});\mathbf{Q}^{t})}.\label{eq:dinkelbach-alpha-parallel}
\end{equation}
It follows from the convergence properties of the Dinkelbach's algorithm
(cf. \cite{Zappone2015}) that
\[
\lim_{\tau\rightarrow\infty}\mathbf{Q}^{\star}(s^{t,\tau})=\mathbb{B}\mathbf{Q}^{t}
\]
at a superlinear convergence rate. Note that $\mathbb{B}\mathbf{Q}^{t}$
is unique, because both $\lim_{\tau\rightarrow\infty}s^{t,\tau}$
and $\mathbf{Q}^{\star}(s^{t,\tau})$ are unique. This iterative procedure
(\ref{eq:dinkelbach-parallel-subproblem})-(\ref{eq:dinkelbach-alpha-parallel})
is nested under Step 1 of Algorithm \ref{alg-GlobalEE} as Steps 1.0-1.3.

In the following, we discuss some properties and implementation aspects
of the proposed Algorithm \ref{alg-GlobalEE}.

The proposed algorithm presents a fast convergence behavior. The approximate
function in (\ref{eq:notation-approximate_function}) is constructed
in the same spirit as \cite{Scutarib,Yang_ConvexApprox} by keeping
as much concavity as possible, namely, $r_{k}(\mathbf{Q}_{k},\mathbf{Q}_{-k})$
in $\mathbf{Q}_{k}$ and $\sum_{j=1}^{K}(P_{0,k}+\rho_{k}\textrm{tr}(\mathbf{Q}_{k}))$
in $\mathbf{Q}$, and linearizing only the nonconcave functions in
the numerator and the nonconvex functions in the denominator, namely,
$\sum_{j\neq k}r_{j}(\mathbf{Q})$ and $\sum_{j=1}^{K}g_{j}(r_{j}(\mathbf{Q}))$.
Besides this, the division operator is also kept. Therefore, the proposed
algorithm is of a best-response nature and expected to exhibit a fast
convergence behavior, as we shall later illustrate numerically.

The proposed algorithm enjoys a low complexity and an easy implementation.
In iterative algorithms, the major computational complexity lies in
solving the approximate problem in each iteration. In the proposed
algorithm, the approximate problem can be decomposed into multiple
independent subproblems and is thus suitable for parallel computation.
The optimal point of each subproblem has a closed-form expression;
by contrast, a generic convex optimization problem must be solved
in each iteration in \cite{Zappone2017,Tervo2017}.

The proposed algorithm presents a broad applicability. Firstly, it
does not require the approximate function to be a global lower bound
of the original function, see, e.g., the sequential programming framework
proposed in \cite{Zappone2017}. Such an approximate function may
not even exist for some choices of the power consumption models. Secondly,
the proposed algorithm is applicable for MIMO systems, where the design
variables are complex-valued matrices, and the rate-dependent processing
power consumption function $g_{k}(\cdot)$ does not have to be convex,
as assumed in \cite{Tervo2017}.

The proposed algorithm can, e.g., be implemented on a central processing
unit which has the channel state information of all direct-link and
cross-link channels, namely, $(\mathbf{H}_{kj})_{j,k}$. In practical
systems, this central unit could be embedded in the centralized radio
access network (CRAN), cf. Figure \ref{fig:System-topology}: each
BS $k$ sends the direct-link channel $\mathbf{H}_{kk}$ and cross-link
channels $(\mathbf{H}_{kj})_{j\neq k}$ to the central unit in the
CRAN. Then the central unit invokes Algorithm \ref{alg-GlobalEE}
and informs each BS $k$ about the optimal transmit covariance matrix
$\mathbf{Q}_{k}$. The incurred latency is mainly due to the signaling
exchange between the central unit and the BSs, and the execution of
the variable updates. Due to the algorithm's low complexity, the central
unit is not required to have a strong computational capability.

\section{\label{sec:Sum-EE}The Proposed Iterative Algorithm for Sum Energy
Efficiency Maximization}

In this section, we propose an iterative algorithm for problem (\ref{eq:original-problem-SumEE}),
which consists in solving a sequence of successively refined approximate
problems. In iteration $t$, we approximate the nonconcave function
$f_{S}(\mathbf{Q})$ with respect to $\mathbf{Q}_{k}$ at the point
$\mathbf{Q}^{t}$ by a function denoted as $\tilde{f}_{S,k}(\mathbf{Q}_{k};\mathbf{Q}^{t})$:
\begin{equation}
\tilde{f}_{S,k}(\mathbf{Q}_{k};\mathbf{Q}^{t})\triangleq\frac{\tilde{r}_{S,k}(\mathbf{Q}_{k};\mathbf{Q}^{t})}{\tilde{p}_{S,k}(\mathbf{Q}_{k};\mathbf{Q}^{t})},\label{eq:sum-ee-approximate-function}
\end{equation}
where \begin{subequations}\label{eq:sum-ee-approximate-function-details}
\begin{align}
\tilde{r}_{S,k}(\mathbf{Q}_{k};\mathbf{Q}^{t}) & \triangleq r_{k}(\mathbf{Q}_{k},\mathbf{Q}_{-k}^{t})+(\mathbf{Q}_{k}-\mathbf{Q}_{k}^{t})\bullet\boldsymbol{\Pi}_{k}(\mathbf{Q}^{t}),\label{eq:sum-ee-approximate-function-numerator-1}\\
\boldsymbol{\Pi}_{k}(\mathbf{Q}^{t}) & \triangleq p_{k}(\mathbf{Q}^{t})\cdot\nabla_{k}\Biggl(\sum_{j\neq k}\frac{r_{j}(\mathbf{Q})}{p_{j}(\mathbf{Q})}\Biggr)\biggr|_{\mathbf{Q}=\mathbf{Q}^{t}},\label{eq:sum-ee-approximate-function-numerator-2}
\end{align}
and
\begin{align}
\tilde{p}_{S,k}(\mathbf{Q}_{k};\mathbf{Q}^{t}) & \triangleq P_{0,k}+\rho_{k}\textrm{tr}(\mathbf{Q}_{k})+g_{k}(r_{k}(\mathbf{Q}^{t}))\nonumber \\
 & \quad\;+(\mathbf{Q}-\mathbf{Q}_{k}^{t})\bullet\nabla_{k}g_{k}(r_{k}(\mathbf{Q}^{t})).\label{eq:sum-ee-approximate-function-denominator}
\end{align}
\end{subequations}In (\ref{eq:sum-ee-approximate-function-numerator-1})-(\ref{eq:sum-ee-approximate-function-numerator-2}),
we fix $\mathbf{Q}_{-k}$ to be $\mathbf{Q}_{-k}=\mathbf{Q}_{-k}^{t}$
in $r_{k}(\mathbf{Q}_{k},\mathbf{Q}_{-k})$ and linearize the nonconcave
function $r_{j}(\mathbf{Q})/p_{j}(\mathbf{Q})$ with respect to $\mathbf{Q}_{k}$.
In (\ref{eq:sum-ee-approximate-function-denominator}), the nonconvex
function $g_{k}(r_{k}(\mathbf{Q}))$ is linearized. As a result, the
numerator and denominator function of $\tilde{f}_{k}(\mathbf{Q}_{k};\mathbf{Q}^{t})$
is concave and convex in $\mathbf{Q}_{k}$, respectively, and $\tilde{f}_{S,k}(\mathbf{Q}_{k};\mathbf{Q}^{t})$
is thus pseudoconcave in $\mathbf{Q}_{k}$. Besides, it is not difficult
to verify that\begin{subequations}\label{eq:SEE-r-equal}
\begin{align}
\tilde{r}_{S,k}(\mathbf{Q}_{k}^{t};\mathbf{Q}^{t}) & =r_{k}(\mathbf{Q}^{t}),\label{eq:SEE-r-equal-value}\\
\left.\nabla_{k}\tilde{r}_{S,k}(\mathbf{Q}_{k};\mathbf{Q}^{t})\right|_{\mathbf{Q}=\mathbf{Q}^{t}} & =\left.\nabla_{k}r_{k}(\mathbf{Q})\right|_{\mathbf{Q}=\mathbf{Q}^{t}}+\boldsymbol{\Pi}_{k}(\mathbf{Q}^{t}),\label{eq:SEE-r-equal-gradient}
\end{align}
\end{subequations}and\begin{subequations}\label{eq:SEE-p-equal}
\begin{align}
\tilde{p}_{S,k}(\mathbf{Q}_{k}^{t};\mathbf{Q}^{t}) & =p_{k}(\mathbf{Q}^{t}),\label{eq:SEE-p-equal-value}\\
\left.\nabla_{k}\tilde{p}_{S,k}(\mathbf{Q}_{k};\mathbf{Q}^{t})\right|_{\mathbf{Q}=\mathbf{Q}^{t}} & =\left.\nabla_{k}p_{k}(\mathbf{Q})\right|_{\mathbf{Q}=\mathbf{Q}^{t}}.\label{eq:SEE-p-equal-gradient}
\end{align}
\end{subequations}Then we can show that $\tilde{f}_{k}(\mathbf{Q}_{k};\mathbf{Q}^{t})$
and $f(\mathbf{Q})$ have the same gradient w.r.t. $\mathbf{Q}_{k}$
at the point $\mathbf{Q}=\mathbf{Q}^{t}$:
\begin{align}
 & \left.\nabla_{k}\tilde{f}_{S,k}(\mathbf{Q}_{k};\mathbf{Q}^{t})\right|_{\mathbf{Q}=\mathbf{Q}^{t}}\nonumber \\
=\; & \frac{\nabla_{k}\tilde{r}_{S,k}(\mathbf{Q}_{k}^{t};\mathbf{Q}^{t})}{\tilde{p}_{S,k}(\mathbf{Q}_{k}^{t};\mathbf{Q}^{t})}-\frac{\tilde{r}_{S,k}(\mathbf{Q}_{k}^{t};\mathbf{Q}^{t})\nabla_{k}\tilde{p}_{S,k}(\mathbf{Q}_{k}^{t};\mathbf{Q}^{t})}{\tilde{p}_{S,k}(\mathbf{Q}_{k}^{t};\mathbf{Q}^{t})^{2}}\nonumber \\
=\; & \frac{\nabla_{k}r_{k}(\mathbf{Q}^{t})+\boldsymbol{\Pi}_{k}(\mathbf{Q}^{t})}{p_{k}(\mathbf{Q}^{t})}-\frac{r_{k}(\mathbf{Q}^{t})\nabla_{k}p_{k}(\mathbf{Q}^{t})}{p_{k}(\mathbf{Q}^{t})^{2}}\nonumber \\
=\; & \frac{\nabla_{k}r_{k}(\mathbf{Q}^{t})}{p_{k}(\mathbf{Q}^{t})}-\frac{r_{k}(\mathbf{Q}^{t})\nabla_{k}p_{k}(\mathbf{Q}^{t})}{p_{k}(\mathbf{Q}^{t})^{2}}\nonumber \\
 & +{\textstyle \sum_{j\neq k}}\nabla_{k}\left(\frac{r_{j}(\mathbf{Q}^{t})}{p_{j}(\mathbf{Q}^{t})}\right)=\left.\nabla_{k}f_{S}(\mathbf{Q})\right|_{\mathbf{Q}=\mathbf{Q}^{t}},\label{eq:SEE-equal-gradient}
\end{align}
where the first and third equality is the expression of $\nabla_{k}\tilde{f}_{S}(\mathbf{Q};\mathbf{Q}^{t})$
and $\nabla_{k}f_{S}(\mathbf{Q})$, respectively, and the second equality
follows from (\ref{eq:SEE-r-equal})-(\ref{eq:SEE-p-equal}).

Given point $\mathbf{Q}^{t}$ in iteration $t$, we define an approximate
problem of the following form:
\begin{align}
\underset{\mathbf{Q}}{\textrm{maximize}}\quad & {\textstyle \sum_{k=1}^{K}}\tilde{f}_{S,k}(\mathbf{Q}_{k};\mathbf{Q}^{t})\nonumber \\
\textrm{subject to}\quad & \mathbf{Q}_{k}\succeq\mathbf{0},\textrm{tr}(\mathbf{Q}_{k})\leq P_{k},\;k=1,\ldots,K,\label{eq:sum-ee-approximate-problem}
\end{align}
and we denote as $\mathbb{B}\mathbf{Q}^{t}=(\mathbb{B}_{k}\mathbf{Q}^{t})_{k=1}^{K}$
the optimal point. Since problem (\ref{eq:sum-ee-approximate-problem})
is well decoupled across different variables, it can be decomposed
into many smaller optimization problems that can be solved in parallel:
\begin{align}
\mathbb{B}\mathbf{Q}^{t} & =\underset{(\mathbf{Q}_{k}\succeq\mathbf{0},\textrm{tr}(\mathbf{Q}_{k})\leq P_{k})_{k=1}^{K}}{\arg\max}\;{\textstyle \sum_{k=1}^{K}}\tilde{f}_{S,k}(\mathbf{Q}_{k};\mathbf{Q}^{t})\label{eq:sum-ee-best-response-1}\\
 & \qquad\qquad\Updownarrow\nonumber \\
\mathbb{B}_{k}\mathbf{Q}^{t} & =\underset{\mathbf{Q}_{k}\succeq\mathbf{0},\textrm{tr}(\mathbf{Q}_{k})\leq P_{k}}{\arg\max}\;\tilde{f}_{S,k}(\mathbf{Q}_{k};\mathbf{Q}^{t}),\;k=1,\ldots,K.\label{eq:sum-ee-best-response}
\end{align}
Note that $\mathbb{B}_{k}\mathbf{Q}^{t}$ is unique, which can be
shown by the same line of argument used in the previous section.

We remark that although $\tilde{f}_{S,k}(\mathbf{Q}_{k};\mathbf{Q}^{t})$
is pseudoconcave in $\mathbf{Q}_{k}$, the approximate function $\sum_{k=1}^{K}\tilde{f}_{S,k}(\mathbf{Q}_{k};\mathbf{Q}^{t})$
in (\ref{eq:sum-ee-approximate-problem}) is not necessarily pseudoconcave
in $\mathbf{Q}$, because, unlike concave functions, the sum of pseudoconcave
functions is not always pseudoconcave. Despite the lack of pseudoconcavity
in the approximate function $\sum_{k=1}^{K}\tilde{f}_{k}(\mathbf{Q}_{k};\mathbf{Q}^{t})$
in (\ref{eq:sum-ee-approximate-problem}), $\mathbb{B}\mathbf{Q}^{t}-\mathbf{Q}^{t}$
is still an ascent direction of the original objective function $f(\mathbf{Q})$
at $\mathbf{Q}=\mathbf{Q}^{t}$. To see this, we note that
\begin{align*}
(\mathbb{B}\mathbf{Q}^{t}-\mathbf{Q}^{t})\bullet\nabla f_{S}(\mathbf{Q}^{t}) & ={\textstyle \sum_{k=1}^{K}}(\mathbb{B}_{k}\mathbf{Q}^{t}-\mathbf{Q}_{k}^{t})\bullet\nabla_{k}f_{S}(\mathbf{Q}^{t})\\
={\textstyle \sum_{k=1}^{K}} & (\mathbb{B}_{k}\mathbf{Q}^{t}-\mathbf{Q}_{k}^{t})\bullet\nabla_{k}\tilde{f}_{S,k}(\mathbf{Q}_{k}^{t};\mathbf{Q}^{t}).
\end{align*}
As both the objective function and the constraint set of the approximate
problem (\ref{eq:sum-ee-approximate-problem}) is well decoupled among
the different block variables and each subproblem (\ref{eq:sum-ee-best-response})
is pseudoconvex, we have $(\mathbb{B}_{k}\mathbf{Q}^{t}-\mathbf{Q}_{k}^{t})\bullet\nabla_{k}\tilde{f}_{S,k}(\mathbf{Q}_{k}^{t};\mathbf{Q}^{t})\geq0$.
This is formally stated in the following proposition.
\begin{prop}
[Stationary point and ascent direction]\label{prop:SEE-ascent-direction}A
point $\mathbf{Q}^{t}$ is a stationary point of (\ref{eq:original-problem-SumEE})
if and only if $\mathbb{B}\mathbf{Q}^{t}=\mathbf{Q}^{t}$. If $\mathbf{Q}^{t}$
is not a stationary point of (\ref{eq:original-problem-SumEE}), then
$\mathbb{B}\mathbf{Q}^{t}-\mathbf{Q}^{t}$ is an ascent direction
of $f_{S}(\mathbf{Q})$ in the sense that
\[
(\mathbb{B}\mathbf{Q}^{t}-\mathbf{Q}^{t})\bullet\nabla f_{S}(\mathbf{Q}^{t})>0.
\]
\end{prop}
\begin{IEEEproof}
It follows from the definition of $\mathbb{B}_{k}\mathbf{Q}^{t}$
in (\ref{eq:sum-ee-best-response}) that
\begin{equation}
\tilde{f}_{k}(\mathbb{B}_{k}\mathbf{Q}^{t};\mathbf{Q}^{t})=\underset{\mathbf{Q}_{k}\succeq\mathbf{0},\textrm{tr}(\mathbf{Q}_{k})\leq P_{k}}{\max}\tilde{f}_{k}(\mathbf{Q}_{k};\mathbf{Q}^{t})\geq\tilde{f}_{k}(\mathbf{Q}_{k}^{t};\mathbf{Q}^{t}),\label{eq:SEE-proof-1}
\end{equation}
and $\mathbb{B}_{k}\mathbf{Q}^{t}=\mathbf{Q}_{k}^{t}$ if equality
holds, as $\mathbb{B}_{k}\mathbf{Q}^{t}$ is unique.

If equality holds in (\ref{eq:SEE-equal-gradient}) for all $k=1,\ldots,K$,
then $\mathbf{Q}^{t}=\mathbb{B}\mathbf{Q}^{t}$ and $\mathbf{Q}^{t}$
is an optimal point of the optimization problem in (\ref{eq:sum-ee-best-response}).
According to the first-order optimality condition, the following inequality
is satisfied:
\[
(\mathbf{Q}_{k}-\mathbf{Q}_{k}^{t})\bullet\nabla_{k}\tilde{f}_{S,k}(\mathbf{Q}_{k}^{t};\mathbf{Q}^{t})\geq0
\]
for all $\mathbf{Q}_{k}$ such that $\mathbf{Q}_{k}\succeq\mathbf{0}$
and $\textrm{tr}(\mathbf{Q}_{k})\leq P_{k}$, $k=1,\ldots,K$. Since
$\nabla_{k}\tilde{f}_{S,k}(\mathbf{Q}_{k}^{t};\mathbf{Q}^{t})=\nabla_{k}f_{S}(\mathbf{Q}^{t})$,
cf. (\ref{eq:SEE-equal-gradient}), the above inequality is
\[
(\mathbf{Q}_{k}-\mathbf{Q}_{k}^{t})\bullet\nabla_{k}f_{S}(\mathbf{Q}^{t})\geq0,\;k=1,\ldots,K.
\]
Adding them up over $k=1,\ldots,K$ yields the first order optimality
condition of the original problem (\ref{eq:original problem-GlobalEE})
and $\mathbf{Q}^{t}$ is thus a stationary point of (\ref{eq:original problem-GlobalEE}).

If strict inequality holds for some $k\in\{1,\ldots,K\}$ in (\ref{eq:SEE-proof-1}),
the pseudoconvexity of the optimization problem in (\ref{eq:sum-ee-best-response})
implies that \begin{subequations}\label{eq:ascent-direction}
\begin{align}
0 & <(\mathbb{B}_{k}\mathbf{Q}^{t}-\mathbf{Q}_{k}^{t})\bullet\nabla_{k}\tilde{f}_{k}(\mathbf{Q}_{k}^{t};\mathbf{Q}^{t})\label{eq:ascnet-direction-approximate}\\
 & =(\mathbb{B}_{k}\mathbf{Q}^{t}-\mathbf{Q}_{k}^{t})\bullet\nabla_{k}f(\mathbf{Q}^{t}),\label{eq:ascent-direction-cartesian}
\end{align}
\end{subequations}where the inequality in (\ref{eq:ascnet-direction-approximate})
comes from the definition of pseudoconcave functions (\ref{eq:def-pseudo-convex-function})
and the equality in (\ref{eq:ascent-direction-cartesian}) comes from
(\ref{eq:SEE-equal-gradient}). Adding up (\ref{eq:ascent-direction})
over all $k=1,\ldots,K$, we obtain
\begin{equation}
(\mathbb{B}\mathbf{Q}^{t}-\mathbf{Q}^{t})\bullet\nabla f_{S}(\mathbf{Q}^{t})>0.\label{eq:sum-EE-ascent-direction}
\end{equation}
The proof is thus completed.
\end{IEEEproof}
According to Proposition \ref{prop:SEE-ascent-direction}, $\mathbb{B}\mathbf{Q}^{t}-\mathbf{Q}^{t}$
is an ascent direction of $f_{S}(\mathbf{Q})$ at $\mathbf{Q}=\mathbf{Q}^{t}$,
and we calculate the stepsize by the successive line search: given
two scalars $0<\alpha<1$ and $0<\beta<1$, $\gamma^{t}$ is set to
be $\gamma^{t}=\beta^{m_{t}}$, where $m_{t}$ is the smallest nonnegative
integer $m$ satisfying the following inequality:
\begin{align}
f_{S}(\mathbf{Q}^{t}+\beta^{m}(\mathbb{B}\mathbf{Q}^{t}-\mathbf{Q}^{t})) & \geq\nonumber \\
f_{S}(\mathbf{Q}^{t})+\alpha\beta^{m} & \nabla f_{S}(\mathbf{Q}^{t})\bullet(\mathbb{B}\mathbf{Q}^{t}-\mathbf{Q}^{t}).\label{eq:line-search-1}
\end{align}
Note that the successive line search is carried out over the original
objective function $f(\mathbf{Q})$ defined in (\ref{eq:original-problem-SumEE}).
After the stepsize $\gamma^{t}$ is found, the variable $\mathbf{Q}$
is updated as
\begin{equation}
\mathbf{Q}^{t+1}=\mathbf{Q}^{t}+\gamma^{t}(\mathbb{B}\mathbf{Q}^{t}-\mathbf{Q}^{t}).\label{eq:sum-EE-variable-update}
\end{equation}

The above steps are formally summarized in Algorithm \ref{alg-SumEE}.
From (\ref{eq:sum-EE-ascent-direction})-(\ref{eq:sum-EE-variable-update})
it can be verified that the sequence $\{f_{S}(\mathbf{Q}^{t})\}_{t}$
is monotonically increasing. Moreover, the sequence $\{\mathbf{Q}^{t}\}$
has a limit point and every limit point is a stationary point of (\ref{eq:original-problem-SumEE}),
whose proof follows the same line of analysis as \cite[Th. 3]{Yang_ConvexApprox}
and thus not duplicated here.

In Step 1 of Algorithm \ref{alg-SumEE}, a constrained pseudoconvex
optimization problem, namely, problem (\ref{eq:sum-ee-best-response}),
must be solved, and we apply the Dinkelbach's algorithm to find $\mathbb{B}_{S,k}\mathbf{Q}^{t}$
iteratively. At iteration $\tau$ of Dinkelbach's algorithm, the following
problem is solved for a given $s_{k}^{t,\tau}$ ($s_{k}^{t,0}$ can
be set to 0):
\begin{align}
\underset{\mathbf{Q}_{k}}{\textrm{maximize}}\quad & \tilde{r}_{S,k}(\mathbf{Q}_{k};\mathbf{Q}^{t})-s_{k}^{t,\tau}\tilde{p}_{S,k}(\mathbf{Q}_{k};\mathbf{Q}^{t})\nonumber \\
\textrm{subject to}\quad & \mathbf{Q}_{k}\succeq\mathbf{0},\,\textrm{tr}(\mathbf{Q}_{k})\leq P_{k},\;\forall k.\label{eq:dinkelbach-parallel-1}
\end{align}
Similar to problem (\ref{eq:dinkelbach-parallel-subproblem}), the
optimal point of problem (\ref{eq:dinkelbach-parallel-1}), denoted
as $\mathbf{Q}_{k}^{\star}(s_{k}^{t,\tau})$, has a closed-form expression
based on the generalized waterfilling solution (cf. (\ref{eq:dinkelbach-parallel-subproblem})
in Section \ref{sec:GlobalEE}). After $(\mathbf{Q}_{k}^{\star}(s_{k}^{t,\tau}))_{k=1}^{K}$
is obtained, $s_{k}^{t,\tau}$ is updated as follows:
\begin{equation}
s_{k}^{t,\tau+1}=\;\frac{\tilde{r}_{S,k}(\mathbf{Q}_{k}^{\star}(s_{k}^{t,\tau});\mathbf{Q}^{t})}{\tilde{p}_{S,k}(\mathbf{Q}^{\star}(s_{k}^{t,\tau});\mathbf{Q}^{t})}.\label{eq:dinkelbach-alpha-parallel-1}
\end{equation}
It follows from the convergence properties of the Dinkelbach's algorithm
that $\lim_{\tau\rightarrow\infty}\mathbf{Q}_{k}^{\star}(s_{k}^{t,\tau})=\mathbb{B}_{k}\mathbf{Q}^{t}$
for all $k$. This iterative procedure (\ref{eq:dinkelbach-parallel-1})-(\ref{eq:dinkelbach-alpha-parallel-1})
is nested under Step 1 of Algorithm \ref{alg-SumEE} as Steps 1.0-1.3.

\begin{algorithm}[t]
\textbf{S0:} $\mathbf{Q}^{0}=\mathbf{0}$, $t=0$, and a stopping
criterion $\varepsilon$.

\textbf{S1: }Compute $\mathbb{B}\mathbf{Q}^{t}$ by solving problem
(\ref{eq:sum-ee-best-response}):

$\qquad$\textbf{S1.0: }$s^{t,0}=0$, $\tau=0$, and a stopping criterion
$\epsilon$.

$\qquad$\textbf{S1.1: }Compute $\mathbf{Q}_{k}^{\star}(s^{t,\tau})$
by (\ref{eq:dinkelbach-parallel-1}).

$\qquad$\textbf{S1.2: }Compute $s^{t,\tau+1}$ by (\ref{eq:dinkelbach-alpha-parallel-1}).

$\qquad$\textbf{S1.3: }If $|s^{t,\tau+1}-s^{t,\tau}|<\epsilon$,
then $\mathbb{B}\mathbf{Q}^{t}=\mathbf{Q}^{\star}(s^{t,\tau})$. $\qquad$$\qquad$$\;$$\;$$\,$Otherwise
$\tau\leftarrow\tau+1$ and go to \textbf{S1.1}.

\textbf{S2: }Compute $\gamma^{t}$ by the successive line search (\ref{eq:line-search-1}).

\textbf{S3: }Update $\mathbf{Q}^{t+1}$ according to (\ref{eq:sum-EE-variable-update}).

\textbf{S4: }If $\left\Vert \mathbb{B}\mathbf{Q}^{t}-\mathbf{Q}^{t}\right\Vert \leq\varepsilon$,
then STOP; otherwise $t\leftarrow t+1$ and $\quad$$\;$$\,$$\,$go
to \textbf{S1}.

\caption{\label{alg-SumEE}The successive pseudoconvex approximation method
for SEE maximization (\ref{eq:original-problem-SumEE})}
\end{algorithm}

The proposed Algorithm \ref{alg-SumEE} for the SEE maximization problem
(\ref{eq:original-problem-SumEE}) has the same attractive features
as those of Algorithm \ref{alg-GlobalEE} for the global EE maximization
problem (\ref{eq:original problem-GlobalEE}), namely, the fast convergence,
the broad applicability and the low complexity; see the discussion
at the end of Sec. \ref{sec:GlobalEE}. We complement the discussion
by emphasizing that Algorithm \ref{alg-SumEE} is the first parallel
best-response Jacobi algorithm designed for the maximization of the
sum EE function, and pseudoconvexity plays a fundamental role that
has not been fully recognized nor exploited by existing techniques.
This also marks a notable relaxation in state-of-the-art convergence
conditions for Jacobi algorithms.

\section{\label{sec:GlobalEE-QoS}The Proposed Algorithm for Global Energy
Efficiency Maximization with QoS Constraints}

In this section, we propose an iterative algorithm to maximize the
GEE subject to the QoS constraints defined in (\ref{eq:original-problem-GlobalEE-w/QoS}).

The nonconcave QoS constraints in (\ref{eq:original-problem-GlobalEE-w/QoS})
make the constraint set nonconvex and Algorithm \ref{alg-GlobalEE}
proposed in Sec. \ref{sec:GlobalEE} for problem (\ref{eq:original problem-GlobalEE})
is no longer applicable, because 1) the approximate problem is difficult
to solve, and 2) the new point updated according to (\ref{eq:variable-update})
is not necessarily feasible. To design an iterative algorithm for
problem (\ref{eq:original-problem-GlobalEE-w/QoS}) that enjoys a
low complexity but at the same time a fast convergence behavior, we
need on the one hand to overcome the nonconcavity/nonconvexity in
the objective function/the constraint set, and, on the other hand,
to preserve the original problem's structure as much as possible.
Towards this end, we extend the successive pseudoconvex approximation
framework developed in \cite{Yang_ConvexApprox} for minimizing a
nonconvex function over a convex constraint set to solve problem (\ref{eq:original-problem-GlobalEE-w/QoS})
where the objective function/the constraint set is nonconcave/nonconvex.

In iteration $t$, the approximate problem defined around the point
$\mathbf{Q}^{t}$ consists of maximizing an approximate function,
denoted as $\tilde{f}(\mathbf{Q};\mathbf{Q}^{t})$, over an approximate
set, denoted as $\tilde{\mathcal{Q}}(\mathbf{Q}^{t})$. We first note
that the nonconcave function $r_{k}(\mathbf{Q})$ in (\ref{eq:original-problem-GlobalEE-w/QoS})
can be rewritten as the difference of two concave functions:
\begin{align*}
r_{k}(\mathbf{Q})=\; & \log\det\left(\mathbf{I}+\mathbf{R}_{k}(\mathbf{Q}_{-k})^{-1}\mathbf{H}_{kk}\mathbf{Q}_{k}\mathbf{H}_{kk}^{H}\right)\\
=\; & \log\det\bigl(\sigma_{k}^{2}\mathbf{I}+{\textstyle \sum_{j=1}^{K}}\mathbf{H}_{kj}\mathbf{Q}_{j}^{H}\mathbf{H}_{kj}^{H}\bigr)\\
 & -\log\det\bigl(\sigma_{k}^{2}\mathbf{I}+{\textstyle \sum_{j\neq k}}\mathbf{H}_{kj}\mathbf{Q}_{j}\mathbf{H}_{kj}^{H}\bigr).
\end{align*}
Introducing auxiliary variables $\mathbf{Y}_{k}$ such that $\mathbf{Y}_{k}={\textstyle \sum_{j=1}^{K}}\mathbf{H}_{kj}\mathbf{Q}_{j}\mathbf{H}_{kj}^{H}$,
we reformulate problem (\ref{eq:original-problem-GlobalEE-w/QoS})
as follows:\begin{subequations}\label{eq:GEE-QoS-new-problem}
\begin{align}
\underset{\mathbf{Q},\mathbf{Y}}{\textrm{maximize}}\quad & f_{G}(\mathbf{Q})\label{eq:GEE-QoS-new-problem-objective}\\
\textrm{subject to}\quad & \mathbf{Q}_{k}\succeq\mathbf{0},\,\textrm{tr}(\mathbf{Q}_{k})\leq P_{k},\label{eq:GEE-QoS-new-problem-constraint1}\\
 & r_{k}^{+}(\mathbf{Y}_{k})-r_{k}^{-}(\mathbf{Q}_{-k})\geq R_{k},\label{eq:GEE-QoS-new-problem-constraint2}\\
 & \mathbf{Y}_{k}={\textstyle \sum_{j=1}^{K}}\mathbf{H}_{kj}\mathbf{Q}_{j}\mathbf{H}_{kj}^{H},\;\forall k,\label{eq:GEE-QoS-new-problem-constraint3}
\end{align}
\end{subequations}where $r_{k}^{+}(\mathbf{Y}_{k})\triangleq\log\det(\sigma_{k}^{2}\mathbf{I}+\mathbf{Y}_{k})$
and $r_{k}^{-}(\mathbf{Q})\triangleq\log\det(\sigma_{k}^{2}\mathbf{I}+{\textstyle \sum_{j\neq k}}\mathbf{H}_{kj}\mathbf{Q}_{j}\mathbf{H}_{kj}^{H})$.
As we will see later, such a reformulation is beneficial because the
resulting approximate problem can be efficiently solved by parallel
algorithms.

\textbf{Approximate function.}\emph{ }The nonconcave numerator function
$\sum_{j=1}^{K}r_{j}(\mathbf{Q})$ is approximated in the same way
as in (\ref{eq:approximate-r}). We also approximate the nonconvex
denominator function $\sum_{j=1}^{K}p_{j}(\mathbf{Q})$ w.r.t. $\mathbf{Q}_{k}$
by $\tilde{p}_{G,k}(\mathbf{Q}_{k};\mathbf{Q}^{t})$ defined in (\ref{eq:approximate-p}).
The approximate function $\tilde{f}(\mathbf{Q};\mathbf{Q}^{t})$ is
of the following form:
\begin{equation}
\tilde{f}_{G}(\mathbf{Q,Y};\mathbf{Q}^{t},\mathbf{Y}^{t})\triangleq\frac{\sum_{k=1}^{K}(\tilde{r}_{G,k}(\mathbf{Q}_{k};\mathbf{Q}^{t})-c\left\Vert \mathbf{Y}_{k}-\mathbf{Y}_{k}^{t}\right\Vert _{F}^{2})}{\sum_{k=1}^{K}\tilde{p}_{G,k}(\mathbf{Q}_{k};\mathbf{Q}^{t})},\label{eq:notation-approximate_function-1}
\end{equation}
with $\mathbf{Y}_{k}^{t}={\textstyle \sum_{j=1}^{K}}\mathbf{H}_{kj}\mathbf{Q}_{j}^{t}\mathbf{H}_{kj}^{H}$,
while $c\geq0$ is a given constant. When $c=0$, the approximate
function (\ref{eq:notation-approximate_function-1}) is the same as
(\ref{eq:notation-approximate_function}). However, when $c>0$, the
quadratic regularization term makes the numerator function strongly
concave in $\mathbf{Y}$ and the benefit will become clear later.
The approximate function has the following important properties:
\begin{itemize}
\item The function $\tilde{f}_{G}(\mathbf{Q,Y};\mathbf{Q}^{t},\mathbf{Y}^{t})$
is pseudoconcave in $(\mathbf{Q,Y})$ for any given and fixed $(\mathbf{Q}^{t},\mathbf{Y}^{t})$.
\item The gradient of $\tilde{f}_{G}(\mathbf{Q,Y};\mathbf{Q}^{t},\mathbf{Y}^{t})$
and that of $f_{G}(\mathbf{Q})$ are identical at the point $(\mathbf{Q}^{t},\mathbf{Y}^{t})$:
\begin{align}
\nabla_{\mathbf{Q}^{\star}}\tilde{f}_{G}(\mathbf{Q},\mathbf{Y};\mathbf{Q}^{t},\mathbf{Y}^{t})\bigr|_{\mathbf{Q}=\mathbf{Q}^{t},\mathbf{Y}=\mathbf{Y}^{t}} & =\nabla_{\mathbf{Q}^{*}}f_{G}(\mathbf{Q}^{t}),\nonumber \\
\nabla_{\mathbf{Y}^{\star}}\tilde{f}_{G}(\mathbf{Q},\mathbf{Y};\mathbf{Q}^{t},\mathbf{Y}^{t})\bigr|_{\mathbf{Q}=\mathbf{Q}^{t},\mathbf{Y}=\mathbf{Y}^{t}}=\mathbf{0} & =\nabla_{\mathbf{Y}^{*}}f_{G}(\mathbf{Q}^{t}).\label{eq:equal-gradient}
\end{align}
\end{itemize}
As we have seen repeatedly, these properties are essential in establishing
the convergence of the proposed algorithm.

\textbf{Approximate set.}\emph{ }It follows from the definition of
concave functions that $r_{k}^{-}(\mathbf{Q})$ is upper bounded by
its first order approximation at the point $\mathbf{Q}^{t}$:
\begin{align}
r_{k}^{-}(\mathbf{Q}) & \leq\,r_{k}^{-}(\mathbf{Q}^{t})\!+\!{\textstyle \sum_{j\neq k}}(\mathbf{Q}_{j}\!-\!\mathbf{Q}_{j}^{t})\!\bullet\!\nabla_{\mathbf{Q}_{j}^{*}}r_{k}^{-}(\mathbf{Q}^{t})\nonumber \\
 & \triangleq\bar{r}_{k}^{-}(\mathbf{Q};\mathbf{Q}^{t}),\label{eq:lower-bound-0}
\end{align}
where
\begin{equation}
r_{k}^{-}(\mathbf{Q}^{t})=\overline{r}_{k}^{-}(\mathbf{Q}^{t};\mathbf{Q}^{t})\textrm{ and }\nabla_{\mathbf{Q}^{*}}r_{k}^{-}(\mathbf{Q}^{t})=\nabla_{\mathbf{Q}^{*}}\overline{r}_{k}^{-}(\mathbf{Q}^{t};\mathbf{Q}^{t}).\label{eq:lower-bound-1}
\end{equation}
Thus $r_{k}^{+}(\mathbf{Y}_{k})-\bar{r}_{k}^{-}(\mathbf{Q};\mathbf{Q}^{t})$
is a global lower bound of $r_{k}(\mathbf{Q})$:
\begin{equation}
r_{k}(\mathbf{Q})=r_{k}^{+}(\mathbf{Y}_{k})-r_{k}^{-}(\mathbf{Q})\geq r_{k}^{+}(\mathbf{Y}_{k})-\bar{r}_{k}^{-}(\mathbf{Q};\mathbf{Q}^{t}),\label{eq:lower-bound-2}
\end{equation}
where equality holds at $\mathbf{Q}=\mathbf{Q}^{t}$.

We then define the (inner) approximate constraint set $\tilde{\mathcal{Q}}(\mathbf{Q}^{t})$
by replacing the nonconcave functions $r_{k}(\mathbf{Q})$ with its
lower bound $\underline{r}_{k}(\mathbf{Q};\mathbf{Q}^{t})$:
\begin{equation}
\tilde{\mathcal{Q}}(\mathbf{Q}^{t})\triangleq\left\{ \negthickspace\negthickspace\begin{array}{cl}
 & \mathbf{Q}_{k}\succeq\mathbf{0},\textrm{tr}(\mathbf{Q}_{k})\leq P_{k},\smallskip\\
(\mathbf{Q,Y}): & r_{k}^{+}(\mathbf{Y}_{k})-\bar{r}_{k}^{-}(\mathbf{Q}_{-k};\mathbf{Q}^{t})\geq R_{k},\\
 & \mathbf{Y}_{k}={\textstyle \sum_{j=1}^{K}}\mathbf{H}_{kj}\mathbf{Q}_{j}\mathbf{H}_{kj}^{H},\forall k
\end{array}\negthickspace\negthickspace\right\} .\label{eq:notation-approximate_feasible_set}
\end{equation}
The set $\tilde{\mathcal{Q}}(\mathbf{Q}^{t})$ is convex as $r_{k}^{+}(\mathbf{Y}_{k})-\bar{r}_{k}(\mathbf{Q};\mathbf{Q}^{t})$
is concave.

\textbf{Approximate problem.}\emph{ }In iteration $t$, the approximate
problem defined at the point $\mathbf{Q}^{t}$ is to maximize the
approximate function $\tilde{f}(\mathbf{Q,Y};\mathbf{Q}^{t},\mathbf{Y}^{t})$
defined in (\ref{eq:notation-approximate_function-1}) over the approximate
set $\tilde{\mathcal{Q}}(\mathbf{Q}^{t})$ defined in and (\ref{eq:notation-approximate_feasible_set}):
\begin{align}
\underset{(\mathbf{Q,Y})\in\tilde{\mathcal{Q}}(\mathbf{Q}^{t})}{\textrm{maximize}}\quad & \tilde{f}_{G}(\mathbf{Q,Y};\mathbf{Q}^{t},\mathbf{Y}^{t}),\label{eq:approximate-problem-2}
\end{align}
and its optimal point is denoted as $(\mathbb{B}_{Q}\mathbf{Q}^{t},\mathbb{B}_{Y}\mathbf{Q}^{t}))$.
Note that its dependence on $\mathbf{Y}^{t}$ is suppressed for notation
simplicity.

It turns out that $\mathbb{B}_{Q}\mathbf{Q}^{t}-\mathbf{Q}^{t}$ is
an ascent direction of the original objective function $f(\mathbf{Q})$
at $\mathbf{Q}=\mathbf{Q}^{t}$, unless $\mathbf{Q}^{t}$ is already
a KKT point\footnote{For an optimization problem with a nonconvex constraint set, a stationary
point is defined as a KKT point, see \cite[Definition 2]{Scutari2017}.} of problem (\ref{eq:original-problem-GlobalEE-w/QoS}), as stated
in the following proposition.
\begin{prop}
[KKT point and ascent direction]\label{prop:ascent-direction-1}A
point $\mathbf{Q}^{t}$ is a KKT point of (\ref{eq:original-problem-GlobalEE-w/QoS})
if and only if $\mathbf{Q}^{t}=\mathbb{B}_{Q}\mathbf{Q}^{t})$. If
$\mathbf{Q}^{t}$ is not a stationary point of (\ref{eq:original-problem-GlobalEE-w/QoS}),
then $\mathbb{B}_{Q}(\mathbf{Q}^{t},\mathbf{Y}^{t})-\mathbf{Q}^{t}$
is an ascent direction of $f_{G}(\mathbf{Q})$ in the sense that
\[
(\mathbb{B}_{Q}\mathbf{Q}^{t}-\mathbf{Q}^{t})\bullet\nabla_{\mathbf{Q}^{*}}f_{G}(\mathbf{Q}^{t})>0.
\]
\end{prop}
\begin{IEEEproof}
See Appendix.
\end{IEEEproof}
Given the ascent direction $\mathbb{B}_{Q}\mathbf{Q}^{t}-\mathbf{Q}^{t}$,
we calculate the stepsize $\gamma^{t}$ by the successive line search
as explained in (\ref{eq:line-search}) and update the variable $\mathbf{Q}$
accordingly. The proposed algorithm is summarized in Algorithm \ref{alg-GEE-w/QoS}
and its convergence properties are given in the following theorem.
\begin{thm}
[Convergence to a KKT point]\label{thm:convergence-1}Given a feasible
initial point $\mathbf{Q}^{0}\in\mathcal{Q}$, the sequence $\{\mathbf{Q}^{t}\}$
generated by Algorithm \ref{alg-GEE-w/QoS} has a limit point, and
every limit point is a KKT point of problem (\ref{eq:original-problem-GlobalEE-w/QoS}).
\end{thm}
\begin{IEEEproof}
Although the constraint set $\mathcal{Q}$ of problem (\ref{eq:original-problem-GlobalEE-w/QoS})
is nonconvex, the sequence $\{\mathbf{Q}^{t}\}$ generated by Algorithm
\ref{alg-GEE-w/QoS} is always feasible. To see this, we check if
$\mathbf{Q}^{t+1}$ satisfies the QoS constraint $r_{k}(\mathbf{Q}^{t+1})\geq R_{k}$:
\begin{align*}
r_{k}(\mathbf{Q}^{t+1})=\; & r_{k}(\mathbf{Q}^{t}+\gamma(\mathbb{B}\mathbf{Q}^{t}-\mathbf{Q}^{t}))\\
\geq\; & r_{k}^{+}(\mathbf{Y}^{t}+\gamma(\mathbb{B}_{Y}\mathbf{Q}^{t}-\mathbf{Y}^{t}))\\
 & -\overline{r}_{k}^{-}(\mathbf{Q}^{t}+\gamma(\mathbb{B}_{Q}\mathbf{Q}^{t}-\mathbf{Q}^{t});\mathbf{Q}^{t})\\
\geq\; & (1-\gamma)(r_{k}^{+}(\mathbf{Y}^{t})-\overline{r}_{k}^{-}(\mathbf{Q}^{t};\mathbf{Q}^{t}))\\
 & +\gamma(r_{k}^{+}(\mathbb{B}_{Y}\mathbf{Q}^{t})-\overline{r}_{k}^{-}(\mathbb{B}_{Q}\mathbf{Q}^{t};\mathbf{Q}^{t}))\\
\geq\; & (1-\gamma)r_{k}(\mathbf{Q}^{t})+\gamma R_{k},
\end{align*}
where the first inequality follows from the fact that $r_{k}^{+}(\mathbf{Y})-\overline{r}_{k}^{-}(\mathbf{Q};\mathbf{Q}^{t})$
is a global lower bound of $r_{k}(\mathbf{Q})$, cf. (\ref{eq:lower-bound-2}),
the second inequality from the concavity of $r_{k}^{+}(\mathbf{Y})-\overline{r}_{k}^{-}(\mathbf{Q};\mathbf{Q}^{t})$,
and the third inequality from the feasibility of $(\mathbb{B}_{Q}\mathbf{Q}^{t},\mathbb{B}_{Y}\mathbf{Q}^{t})$,
i.e., $(\mathbb{B}_{Q}\mathbf{Q}^{t},\mathbb{B}_{Y}\mathbf{Q}^{t})\in\tilde{\mathcal{Q}}(\mathbf{Q}^{t})$.
Therefore $r_{k}(\mathbf{Q}^{t+1})\geq R_{k}$ if $r_{k}(\mathbf{Q}^{t})\geq R_{k}$.
Since $\mathbf{Q}^{0}$ is feasible, $\mathbf{Q}^{t+1}$ is feasible
by induction.

Since the constraint set $\mathcal{Q}$ is closed and bounded, the
sequence $\{\mathbf{Q}^{t}\}_{t}$ is bounded and thus has a limit
point. The proof for the latter argument follows the same line of
analysis as \cite[Theorem 1]{Yang_ConvexApprox}.
\end{IEEEproof}
\textbf{On solving the approximate problem (\ref{eq:approximate-problem-2}).}
Proposition \ref{prop:ascent-direction-1} and Theorem \ref{thm:convergence-1}
hold for any choice of nonnegative $c$, even when $c=0$. Since problem
(\ref{eq:approximate-problem-2}) is pseudoconcave, its globally optimal
point can be found either by standard gradient-based methods or by
the interior-point method proposed in \cite{Freund1994}.

The choice of a positive $c$ brings numerical benefits when we apply
the Dinkelbach's algorithm to solve problem (\ref{eq:approximate-problem-2})
iteratively. At iteration $\tau$ of Dinkelbach's algorithm, the following
problem is solved for a given and fixed $s^{t,\tau}$:\begin{subequations}\label{eq:GEE-QoS-dinkelbach-parallel}
\begin{align}
\underset{\mathbf{Q},\mathbf{Y}}{\textrm{maximize}}\quad & {\textstyle \sum_{k=1}^{K}}\Bigl(\tilde{r}_{G,k}(\mathbf{Q}_{k};\mathbf{Q}^{t})-c\left\Vert \mathbf{Y}_{k}-\mathbf{Y}_{k}^{t}\right\Vert _{F}^{2}\Bigr)\nonumber \\
 & -s^{t,\tau}{\textstyle \sum_{k=1}^{K}}\tilde{p}_{G,k}(\mathbf{Q}_{k};\mathbf{Q}^{t})\label{eq:GEE-QoS-dinkelbach-parallel-1}\\
\textrm{subject to}\quad & \mathbf{Q}_{k}\succeq\mathbf{0},\,\textrm{tr}(\mathbf{Q}_{k})\leq P_{k},\label{eq:GEE-QoS-dinkelbach-parallel-2}\\
 & r_{k}^{+}(\mathbf{Y}_{k})-\bar{r}_{k}^{-}(\mathbf{Q}_{-k};\mathbf{Q}^{t})\geq R_{k},\label{eq:GEE-QoS-dinkelbach-parallel-3}\\
 & \mathbf{Y}_{k}={\textstyle \sum_{j=1}^{K}}\mathbf{H}_{kj}\mathbf{Q}_{j}\mathbf{H}_{kj}^{H},\,k=1,\ldots,K.\label{eq:GEE-QoS-dinkelbach-parallel-4}
\end{align}
\end{subequations}We denote the solution of problem (\ref{eq:GEE-QoS-dinkelbach-parallel})
as $(\mathbf{Q}^{\star}(s^{t,\tau}),\mathbf{Y}^{\star}(s^{t,\tau})$.
Then $s^{t,\tau}$ is updated as follows:
\begin{equation}
s^{t,\tau+1}=\frac{\sum_{k=1}^{K}\tilde{r}_{G,k}(\mathbf{Q}_{k}^{\star}(s^{t,\tau});\mathbf{Q}^{t})-c\left\Vert \mathbf{Y}_{k}^{\star}(s^{t,\tau})-\mathbf{Y}_{k}^{t}\right\Vert _{F}^{2}}{\sum_{k=1}^{K}\tilde{p}_{G,k}(\mathbf{Q}_{k}^{\star}(s^{t,\tau});\mathbf{Q}^{t})}.\label{eq:dinkelbach-alpha-parallel-2}
\end{equation}
It follows from the convergence properties of the Dinkelbach's algorithm
that $\lim_{\tau\rightarrow\infty}\mathbf{Q}^{\star}(s^{t,\tau})=\mathbb{B}_{Q}\mathbf{Q}^{t}$
and $\lim_{\tau\rightarrow\infty}\mathbf{Y}^{\star}(s^{t,\tau})=\mathbb{B}_{Y}\mathbf{Q}^{t}$.
This iterative procedure (\ref{eq:GEE-QoS-dinkelbach-parallel})-(\ref{eq:dinkelbach-alpha-parallel-2})
is nested under Step 1 of Algorithm \ref{alg-GEE-w/QoS}.

\begin{algorithm}[tbh]
\textbf{S0:} $\mathbf{Q}^{0}\in\mathcal{Q}$, $t=0$, and a stopping
criterion $\varepsilon$.

\textbf{S1: }Compute $\mathbb{B}_{Q}\mathbf{Q}^{t}$ by solving problem
(\ref{eq:approximate-problem-2}):

$\qquad$\textbf{S1.0: }$s^{t,0}=0$, $\tau=0$, and a stopping criterion
$\epsilon$.

$\qquad$\textbf{S1.1: }Compute $\mathbf{Q}^{\star}(s^{t,\tau})$
by solving problem (\ref{eq:GEE-QoS-dinkelbach-parallel}):

$\qquad$$\qquad$\textbf{S1.1.0:} $\upsilon=0$, $\boldsymbol{\lambda}=\mathbf{0}$,
$\boldsymbol{\Sigma}=\mathbf{0}$, and a stopping criterion $\sigma$.

$\qquad$$\qquad$\textbf{S1.1.1:} Compute $\mathbf{Q}_{k}^{L}(\lambda_{k}^{\upsilon})$
and $\mathbf{Y}_{k}^{L}(\boldsymbol{\Sigma}_{k}^{\upsilon})$ by (\ref{eq:dual-1})
and (\ref{eq:dual-2}), respectively, for all $k=1,\ldots,K$.

$\qquad$$\qquad$\textbf{S1.1.2:} Update $\lambda_{k}$ and $\boldsymbol{\Sigma}_{k}$
by (\ref{eq:dual-variable-update}) for all $k$.

$\qquad$$\qquad$\textbf{S1.1.3:} If $\left\Vert (\boldsymbol{\lambda}^{\upsilon+1},\boldsymbol{\Sigma}^{\upsilon+1})-(\boldsymbol{\lambda}^{\upsilon},\boldsymbol{\Sigma}^{\upsilon})\right\Vert \leq\sigma$,
then $\mathbf{Q}_{k}^{\star}(s^{t,\tau})=\mathbf{Q}_{k}^{L}(\lambda_{k}^{\upsilon})$.
Otherwise $\upsilon\leftarrow\upsilon+1$ and go to \textbf{S1.1.1}.

$\qquad$\textbf{S1.2: }Compute $s^{t,\tau+1}$ by (\ref{eq:dinkelbach-alpha-parallel-2}).

$\qquad$\textbf{S1.3: }If $|s^{t,\tau+1}-s^{t,\tau}|<\epsilon$,
then $\mathbb{B}_{Q}\mathbf{Q}^{t}=\mathbf{Q}^{\star}(s^{t,\tau})$.
$\qquad$$\qquad$$\;$$\;$$\,$Otherwise $\tau\leftarrow\tau+1$
and go to \textbf{S1.1}.

\textbf{S2: }Compute $\gamma^{t}$ by the successive line search (\ref{eq:line-search}).

\textbf{S3: }Update $\mathbf{Q}$ and $\mathbf{Y}$ by $\mathbf{Q}^{t+1}=\mathbf{Q}^{t}+\gamma^{t}(\mathbb{B}_{Q}\mathbf{Q}^{t}-\mathbf{Q}^{t})$
and $\mathbf{Y}^{t+1}=\mathbf{Y}^{t}+\gamma^{t}(\mathbb{B}_{Y}\mathbf{Q}^{t}-\mathbf{Y}^{t})$,
respectively.

\textbf{S3: }If $\left\Vert \mathbb{B}\mathbf{Q}^{t}-\mathbf{Q}^{t}\right\Vert \leq\varepsilon$,
then STOP; otherwise$t\leftarrow t+1$ and go to \textbf{S1}.

\caption{\label{alg-GEE-w/QoS}The successive pseudoconvex approximation method
for GEE maximization with QoS constraints (\ref{eq:original-problem-GlobalEE-w/QoS})}
\end{algorithm}

\textbf{On solving problem (\ref{eq:GEE-QoS-dinkelbach-parallel})}\textbf{\emph{.}}
Problem (\ref{eq:GEE-QoS-dinkelbach-parallel}) is convex and the
coupling constraints have a separable structure, which can readily
be exploited in the standard dual decomposition method. To see this,
the Lagrangian of (\ref{eq:GEE-QoS-dinkelbach-parallel}) is:
\begin{align}
L(\mathbf{Q},\mathbf{Y},\boldsymbol{\lambda}\, & ,\boldsymbol{\Sigma};\mathbf{Q}^{t},\mathbf{Y}^{t},s^{t,\tau})={\textstyle \sum_{k=1}^{K}}\tilde{r}_{G,k}(\mathbf{Q}_{k};\mathbf{Q}^{t})\nonumber \\
 & -{\textstyle \sum_{k=1}^{K}}\bigl(c\bigl\Vert\mathbf{Y}_{k}-\mathbf{Y}_{k}^{t}\bigr\Vert_{F}^{2}+s^{t,\tau}\tilde{p}_{G,k}(\mathbf{Q}_{k};\mathbf{Q}^{t})\bigr)\nonumber \\
 & -{\textstyle \sum_{k=1}^{K}}\boldsymbol{\Sigma}_{k}\bullet\bigl(\mathbf{Y}_{k}-{\textstyle \sum_{j=1}^{K}}\mathbf{H}_{kj}\mathbf{Q}_{j}\mathbf{H}_{kj}^{H}\bigr)\nonumber \\
 & +{\textstyle \sum_{k=1}^{K}}\lambda_{k}(r_{k}^{+}(\mathbf{Y}_{k})-\bar{r}_{k}^{-}(\mathbf{Q}_{-k};\mathbf{Q}^{t})-R_{k}),\label{eq:Lagrangian-1}
\end{align}
where $\lambda_{k}$ and $\boldsymbol{\Sigma}_{k}$ are the Lagrange
multipliers associated with the constraints (\ref{eq:GEE-QoS-dinkelbach-parallel-3})-(\ref{eq:GEE-QoS-dinkelbach-parallel-4}).
The dual function $d(\boldsymbol{\lambda},\boldsymbol{\Sigma})$ is
\begin{equation}
d(\boldsymbol{\lambda},\boldsymbol{\Sigma})=\underset{(\mathbf{Q}_{k}\succeq\mathbf{0},\textrm{tr}(\mathbf{Q}_{k})\leq P_{k},\mathbf{Y}_{k}\succeq\mathbf{0})_{k=1}^{K}}{\max}L(\mathbf{Q},\mathbf{Y},\boldsymbol{\lambda},\boldsymbol{\Sigma}),\label{eq:dual-function-1}
\end{equation}
where the dependence of $L(\mathbf{Q,Y},\boldsymbol{\lambda,\Sigma})$
on $(\mathbf{Q}^{t},\mathbf{Y}^{t},s^{t,\tau})$ is dropped in (\ref{eq:dual-function-1})
for notation simplicity. The dual problem of (\ref{eq:GEE-QoS-dinkelbach-parallel})
is
\begin{equation}
\underset{\boldsymbol{\lambda}\geq\mathbf{0},\boldsymbol{\Sigma}}{\textrm{minimize}}\quad d(\boldsymbol{\lambda},\boldsymbol{\Sigma}).\label{eq:dual-problem-1}
\end{equation}
Since the Lagrangian $L(\mathbf{Q},\mathbf{Y},\boldsymbol{\lambda},\boldsymbol{\Sigma})$
is well decoupled across different variables for fixed dual variable
$(\boldsymbol{\lambda},\boldsymbol{\Sigma})$, the maximization problem
in (\ref{eq:dual-function-1}) can be decomposed into many smaller
optimization problems that can be solved in parallel: for all $k=1,\ldots,K$,
\begin{align}
\mathbf{Q}_{k}^{L}(\lambda_{k})\triangleq\quad & \underset{\mathbf{Q}_{k}\succeq\mathbf{0},\textrm{tr}(\mathbf{Q}_{k})\leq P_{k}}{\arg\max}\nonumber \\
\Bigl\{\tilde{r}_{G,k}(\mathbf{Q}_{k};\mathbf{Q}^{t}) & -s^{t,\tau}\tilde{p}_{G,k}(\mathbf{Q}_{k};\mathbf{Q}^{t})-\overline{r}_{k}^{-}(\mathbf{Q}_{-k};\mathbf{Q}^{t})\Bigr\}\label{eq:dual-1}
\end{align}
and
\begin{equation}
\mathbf{Y}^{L}(\boldsymbol{\Sigma}_{k})\triangleq\underset{\mathbf{Y}_{k}\succeq\mathbf{0}}{\arg\max}\left\{ \begin{array}{l}
\lambda_{k}\log\det(\sigma_{k}^{2}\mathbf{I}+\mathbf{Y}_{k})-\boldsymbol{\Sigma}_{k}\bullet\mathbf{Y}_{k}\smallskip\\
-c\left\Vert \mathbf{Y}_{k}-\mathbf{Y}_{k}^{t}\right\Vert _{F}^{2}
\end{array}\right\} ,\label{eq:dual-2}
\end{equation}
where ``L'' in the superscript stands for ``Lagrangian''. Since
$c>0$, $\mathbf{Q}_{k}^{L}(\lambda_{k})$ in (\ref{eq:dual-1}) and
$\mathbf{Y}^{L}(\boldsymbol{\Sigma}_{k})$ in (\ref{eq:dual-2}) exist
and are unique, and they have a closed-form expression, cf. \cite[Lem. 2]{Yang_stochastic}
and \cite[Lem. 7]{Scutari_NonconvexApplications}.

The dual problem (\ref{eq:dual-problem-1}) can be solved by the gradient
projection algorithm and its gradient of $d(\boldsymbol{\lambda},\boldsymbol{\Sigma})$
is
\begin{align*}
\nabla_{\lambda_{k}}d(\boldsymbol{\lambda},\boldsymbol{\Sigma}) & =\log\det(\sigma_{k}^{2}\mathbf{I}\!+\!\mathbf{Y}_{k}^{L}(\boldsymbol{\Sigma}_{k}))\!-\!\bar{r}_{k}^{-}(\mathbf{Q}_{-k}^{L}(\lambda_{k}))\!-\!R_{k},\\
\nabla_{\boldsymbol{\Sigma}^{*}}d(\boldsymbol{\lambda},\boldsymbol{\Sigma}) & ={\textstyle \sum_{j=1}^{K}}\mathbf{H}_{kj}\mathbf{Q}_{j}^{L}(\lambda_{k})\mathbf{H}_{kj}^{H}-\mathbf{Y}_{k}^{L}(\boldsymbol{\Sigma}_{k}).
\end{align*}
In iteration $\upsilon$ to solve problem (\ref{eq:dual-problem-1}),
the dual variable is updated as follows:\begin{subequations}\label{eq:dual-variable-update}
\begin{align}
\lambda_{k}^{t,\tau,\upsilon+1} & =\Bigl[\lambda_{k}^{t,\tau,\upsilon}+\zeta^{t,\tau,\upsilon}\nabla_{\lambda_{k}}d(\boldsymbol{\lambda}^{t,\tau,\upsilon},\boldsymbol{\Sigma}^{t,\tau,\upsilon})\Bigr]^{+},\label{eq:dual-variable-update-1}\\
\boldsymbol{\Sigma}_{k}^{t,\tau,\upsilon+1} & =\boldsymbol{\Sigma}_{k}^{t,\tau,\upsilon}+\zeta^{t,\tau,\upsilon}\nabla_{\boldsymbol{\Sigma}^{*}}d(\boldsymbol{\lambda}^{t,\tau,\upsilon},\boldsymbol{\Sigma}^{t,\tau,\upsilon}),\label{eq:dual-variable-update-2}
\end{align}
\end{subequations}where $\boldsymbol{\lambda}^{t,\tau,0}$ and $\boldsymbol{\Sigma}^{t,\tau,0}$
can be set to $\mathbf{0}$. If $c>0$ and the stepsizes $\{\zeta^{t,\tau,\upsilon}\}_{\upsilon}$
are properly selected, e.g., $\sum_{\upsilon}\zeta^{t,\tau,\upsilon}=\infty$
and $\sum_{\upsilon}(\zeta^{t,\tau,\upsilon})^{2}<\infty$, then $\lim_{\upsilon\rightarrow\infty}\mathbf{Q}^{l}(\lambda_{k}^{t,\tau,\upsilon})=\mathbf{Q}^{d}(\alpha^{t,\tau})$
and $\lim_{\upsilon\rightarrow\infty}\mathbf{Y}^{l}(\lambda_{k}^{t,\tau,\upsilon})=\mathbf{Y}^{d}(\alpha^{t,\tau})$.
This iterative procedure (\ref{eq:dual-1})-(\ref{eq:dual-variable-update})
is nested under Step 1.1 of Algorithm \ref{alg-GEE-w/QoS}.

The Algorithm \ref{alg-GEE-w/QoS} consists of three layers: the outer
layer with index $t$, middle layer with index $\tau$, and inner
layer with index $\upsilon$. The relationship of different layers
is given as follows: $\mathbf{Q}^{\star}=\lim_{t\rightarrow\infty}\lim_{\tau\rightarrow\infty}\lim_{\upsilon\rightarrow\infty}\mathbf{Q}^{l}(\lambda_{k}^{t,\tau,\upsilon})$,
where $\mathbf{Q}^{\star}$ is a KKT point of (\ref{eq:original-problem-GlobalEE-w/QoS})
and the limit with respect to $t$ is in the sense of subsequence
convergence specified by Theorem \ref{thm:convergence-1}. Note that
although the proposed algorithm consists of three layers, its convergence
speed is not negatively affected, because all updates have closed-form
expressions and both the middle and inner layers converge very fast.
Typically convergence is observed after a few iterations.

\section{\label{sec:SEE-w/QoS}The Proposed Algorithm for Sum Energy Efficiency
Maximization with QoS Constraints}

In this section, we propose an iterative algorithm to maximize the
SEE subject to the QoS constraints defined in (\ref{eq:original-problem-SumEE-w/QoS}).

Given $\mathbf{Q}^{t}$ at iteration $t$, it is tempting to define
the approximate function as (\ref{eq:sum-ee-approximate-problem}),
which is proposed for problem (\ref{eq:original-problem-SumEE}),
where the approximate function is the sum of multiple component functions
$(\tilde{f}_{S,k}(\mathbf{Q}_{k};\mathbf{Q}^{t}))_{k=1}^{K}$, while
each component function $\tilde{f}_{S,k}(\mathbf{Q}_{k};\mathbf{Q}^{t})$
is pseudoconcave in $\mathbf{Q}_{k}$. However, different from problem
(\ref{eq:original-problem-SumEE}), the QoS constraints introduce
coupling among different optimization variables $(\mathbf{Q}_{k})_{k=1}^{K}$
in the constraint set in problem (\ref{eq:original-problem-SumEE-w/QoS}),
making it impossible to decompose the approximate problem into multiple
independent pseudoconvex optimization problems, cf. (\ref{eq:sum-ee-best-response-1})-(\ref{eq:sum-ee-best-response}).
To overcome this difficulty, we define an approximate function that
is concave as concavity is preserved under addition and a concave
function is also pseudoconcave.

Firstly, we reformulate problem (\ref{eq:original-problem-SumEE-w/QoS})
as follows:\begin{subequations}\label{eq:SEE-QoS-new-problem}
\begin{align}
\underset{\mathbf{Q},\mathbf{Y}}{\textrm{maximize}}\quad & f_{S}(\mathbf{Q})\label{eq:SEE-QoS-new-problem-objective}\\
\textrm{subject to}\quad & \mathbf{Q}_{k}\succeq\mathbf{0},\,\textrm{tr}(\mathbf{Q}_{k})\leq P_{k},\label{eq:SEE-QoS-new-problem-constraint1}\\
 & r_{k}^{+}(\mathbf{Y}_{k})-r_{k}^{-}(\mathbf{Q}_{-k})\geq R_{k},\label{eq:SEE-QoS-new-problem-constraint2}\\
 & \mathbf{Y}_{k}={\textstyle \sum_{j=1}^{K}}\mathbf{H}_{kj}\mathbf{Q}_{j}\mathbf{H}_{kj}^{H},\;\forall k,\label{eq:SEE-QoS-new-problem-constraint3}
\end{align}
\end{subequations}On the one hand, we approximate the original objective
function $f_{S}(\mathbf{Q})$ by an approximate function $\tilde{f}_{S}(\mathbf{Q},\mathbf{Y};\mathbf{Q}^{t})$:\begin{subequations}\label{eq:SEE-wQoS-approximate-function}
\begin{align}
\tilde{f}_{S}(\mathbf{Q},\mathbf{Y};\mathbf{Q}^{t})=\; & {\textstyle \sum_{k=1}^{K}}\tilde{f}_{S,k}(\mathbf{Q}_{k},\mathbf{Y}_{k};\mathbf{Q}^{t}),\label{eq:SEE-wQoS-approximate-function-1}\\
\tilde{f}_{S,k}(\mathbf{Q}_{k},\mathbf{Y}_{k};\mathbf{Q}^{t})\triangleq\; & \frac{r_{k}(\mathbf{Q}_{k},\mathbf{Q}_{-k}^{t})}{p_{k}(\mathbf{Q}^{t})}+(\mathbf{Q}_{k}-\mathbf{Q}_{k}^{t})\bullet\boldsymbol{\Pi}_{k}(\mathbf{Q}^{t})\nonumber \\
 & -c\left\Vert \mathbf{Y}_{k}-\mathbf{Y}_{k}^{t}\right\Vert _{F}^{2},\label{eq:SEE-wQoS-approximate-function-2}
\end{align}
\end{subequations}where
\[
\boldsymbol{\Pi}_{k}(\mathbf{Q})=-\frac{r_{k}(\mathbf{Q}^{t})}{p_{k}(\mathbf{Q}^{t})^{2}}\nabla_{\mathbf{Q}_{k}^{*}}p_{k}(\mathbf{Q}^{t})+{\textstyle \sum_{j\neq k}}\nabla_{\mathbf{Q}_{k}^{*}}\left(\frac{r_{j}(\mathbf{Q}^{t})}{p_{j}(\mathbf{Q}^{t})}\right).
\]
In contrast to (\ref{eq:sum-ee-approximate-function}), $\tilde{f}_{S,k}(\mathbf{Q}_{k},\mathbf{Y}_{k};\mathbf{Q}^{t})$
defined in (\ref{eq:SEE-wQoS-approximate-function-2}) is no longer
a fractional function, and it is concave in $(\mathbf{Q}_{k},\mathbf{Y}_{k})$.
Therefore, $\tilde{f}_{S}(\mathbf{Q},\mathbf{Y};\mathbf{Q}^{t})$
is concave in $(\mathbf{Q,Y})$. Furthermore, its gradient at the
point $(\mathbf{Q,Y})=(\mathbf{Q}^{t},\mathbf{Y}^{t})$ is the same
as that of the original function $f_{S}(\mathbf{Q})$:
\begin{align*}
 & \left.\nabla_{\mathbf{Q}_{k}^{*}}\tilde{f}_{S}(\mathbf{Q,Y};\mathbf{Q}^{t})\right|_{\mathbf{Q}=\mathbf{Q}^{t}}=\left.\nabla_{\mathbf{Q}_{k}^{*}}\tilde{f}_{S,k}(\mathbf{Q}_{k},\mathbf{Y}_{k};\mathbf{Q}^{t})\right|_{\mathbf{Q}=\mathbf{Q}^{t}}\\
=\; & \frac{\nabla_{\mathbf{Q}_{k}^{*}}r_{k}(\mathbf{Q}^{t})}{p_{k}(\mathbf{Q}^{t})}-\frac{r_{k}(\mathbf{Q}^{t})}{p_{k}(\mathbf{Q}^{t})^{2}}\nabla_{\mathbf{Q}_{k}^{*}}p_{k}(\mathbf{Q}^{t})+\sum_{j\neq k}\nabla_{\mathbf{Q}_{k}^{*}}\frac{r_{j}(\mathbf{Q}^{t})}{p_{j}(\mathbf{Q}^{t})}\\
=\; & \left.\nabla_{\mathbf{Q}_{k}^{*}}\tilde{f}_{S}(\mathbf{Q})\right|_{\mathbf{Q}=\mathbf{Q}^{t}},
\end{align*}
and $\left.\nabla_{\mathbf{Y}_{k}^{*}}\tilde{f}_{S}(\mathbf{Q,Y};\mathbf{Q}^{t})\right|_{\mathbf{Y}=\mathbf{Y}^{t}}=\mathbf{0}=\nabla_{\mathbf{Y}_{k}^{*}}\tilde{f}_{S}(\mathbf{Q})$.

\begin{algorithm}[t]
\textbf{S0:} $\mathbf{Q}^{0}=\mathbf{0}$, $t=0$, and a stopping
criterion $\varepsilon$.

\textbf{S1: }Compute $(\mathbb{B}_{Q}\mathbf{Q}^{t},\mathbb{B}_{Y}\mathbf{Q}^{t})$
by problem (\ref{eq:SEE-wQoS-approximate-problem}):

\textbf{S2: }Compute $\gamma^{t}$ by the successive line search (\ref{eq:line-search-1}).

\textbf{S3: }Update $\mathbf{Q}$ and $\mathbf{Y}$ by $\mathbf{Q}^{t+1}=\mathbf{Q}^{t}+\gamma^{t}(\mathbb{B}_{Q}\mathbf{Q}^{t}-\mathbf{Q}^{t})$
and $\mathbf{Y}^{t+1}=\mathbf{Y}^{t}+\gamma^{t}(\mathbb{B}_{Y}\mathbf{Q}^{t}-\mathbf{Y}^{t})$,
respectively.

\textbf{S4: }If $\left\Vert \mathbb{B}_{Q}\mathbf{Q}^{t}-\mathbf{Q}^{t}\right\Vert \leq\varepsilon$,
then STOP; otherwise $t\leftarrow t+1$ and go to \textbf{S1}.

\caption{\label{alg-SEE-w/QoS}The successive pseudoconvex approximation method
for SEE maximization with QoS constraints (\ref{eq:original-problem-SumEE-w/QoS})}
\end{algorithm}

On the other hand, the nonconvex constraint set in (\ref{eq:SEE-QoS-new-problem})
is approximated by its inner approximation $\tilde{\mathcal{Q}}(\mathbf{Q}^{t})$
defined in (\ref{eq:notation-approximate_feasible_set}). Then in
iteration $t$, the approximate problem consists of maximizing the
approximate function $\tilde{f}_{S}(\mathbf{Q,Y};\mathbf{Q}^{t})$
over the approximate set $\tilde{\mathcal{Q}}(\mathbf{Q}^{t})$:
\begin{align}
\underset{\mathbf{Q,Y}}{\textrm{maximize}}\quad & {\textstyle \sum_{k=1}^{K}}\tilde{f}_{S,k}(\mathbf{Q}_{k};\mathbf{Q}^{t})\nonumber \\
\textrm{subject to}\quad & \mathbf{Q}_{k}\succeq\mathbf{0},\textrm{tr}(\mathbf{Q}_{k})\leq P_{k},\nonumber \\
 & r_{k}^{+}(\mathbf{Y}_{k})-\overline{r}_{k}^{-}(\mathbf{Q}_{k};\mathbf{Q}^{t})\geq R_{k},\nonumber \\
 & \mathbf{Y}_{k}={\textstyle \sum_{j=1}^{K}}\mathbf{H}_{kj}\mathbf{Q}_{j}\mathbf{H}_{kj}^{H},\;k=1,\ldots,K,\label{eq:SEE-wQoS-approximate-problem}
\end{align}
and let $(\mathbb{B}_{Q}\mathbf{Q}^{t},\mathbb{B}_{Y}\mathbf{Q}^{t})$
denote the optimal point. By following the same line of analysis in
Proposition \ref{prop:ascent-direction-1}, we can show that $\mathbb{B}_{Q}\mathbf{Q}^{t}-\mathbf{Q}^{t}$
is an ascent direction of $f_{S}(\mathbf{Q})$ at $\mathbf{Q}=\mathbf{Q}^{t}$,
unless $\mathbf{Q}^{t}$ is already a KKT point of (\ref{eq:SEE-QoS-new-problem}).
To update the variable, the stepsize could be calculated by the successive
line search as explained in (\ref{eq:line-search-1}). Following the
same line of analysis in Theorem \ref{thm:convergence-1}, we could
claim that the sequence $\{\mathbf{Q}^{t}\}$ has a limit point and
any limit point is a KKT point of (\ref{eq:SEE-QoS-new-problem}).

The above iterative procedure is summarized in Algorithm \ref{alg-SEE-w/QoS}.
In Step 1, the convex optimization problem (\ref{eq:SEE-wQoS-approximate-problem})
is solved. As its objective function and constraint set have a separable
structure, (\ref{eq:SEE-wQoS-approximate-problem}) could be solved
by parallel algorithms based on the dual decomposition. The discussion
is similar to that of problem (\ref{eq:GEE-QoS-dinkelbach-parallel})
and thus omitted here.

\section{\label{sec:Simulations}Simulations}

In this section, we compare numerically the proposed algorithms with
state-of-the-art algorithms. In particular, we consider a 7-user MIMO
IC, where the number of transmit antennas is $M_{T}=8$ and the number
of receive antennas is $M_{R}=4$. The power dissipated in hardware
is $P_{0,k}=10$dB, and the power budget normalized by the number
of transmit antennas is 10dB, i.e., $P_{k}/M_{T}=$10dB. The inverse
of the power amplifier efficiency is $\rho=2.6$, the noise covariance
is $\sigma^{2}=1$, the antenna gain is 16dB, and the path loss exponent
is 2. The results are averaged over 20 i.i.d. random channel realizations.
All algorithms are tested under identical conditions under Matlab
R2017a on a PC equipped with an operating system of Windows 10 64-bit,
an Intel i7-7600U 2.80GHz CPU, and a 16GB RAM. All of the Matlab codes
are available online at \url{https://wwwen.uni.lu/snt/people/yang_yang}.

\subsection{GEE Maximization}

\begin{figure}[t]
\center\includegraphics[scale=0.6]{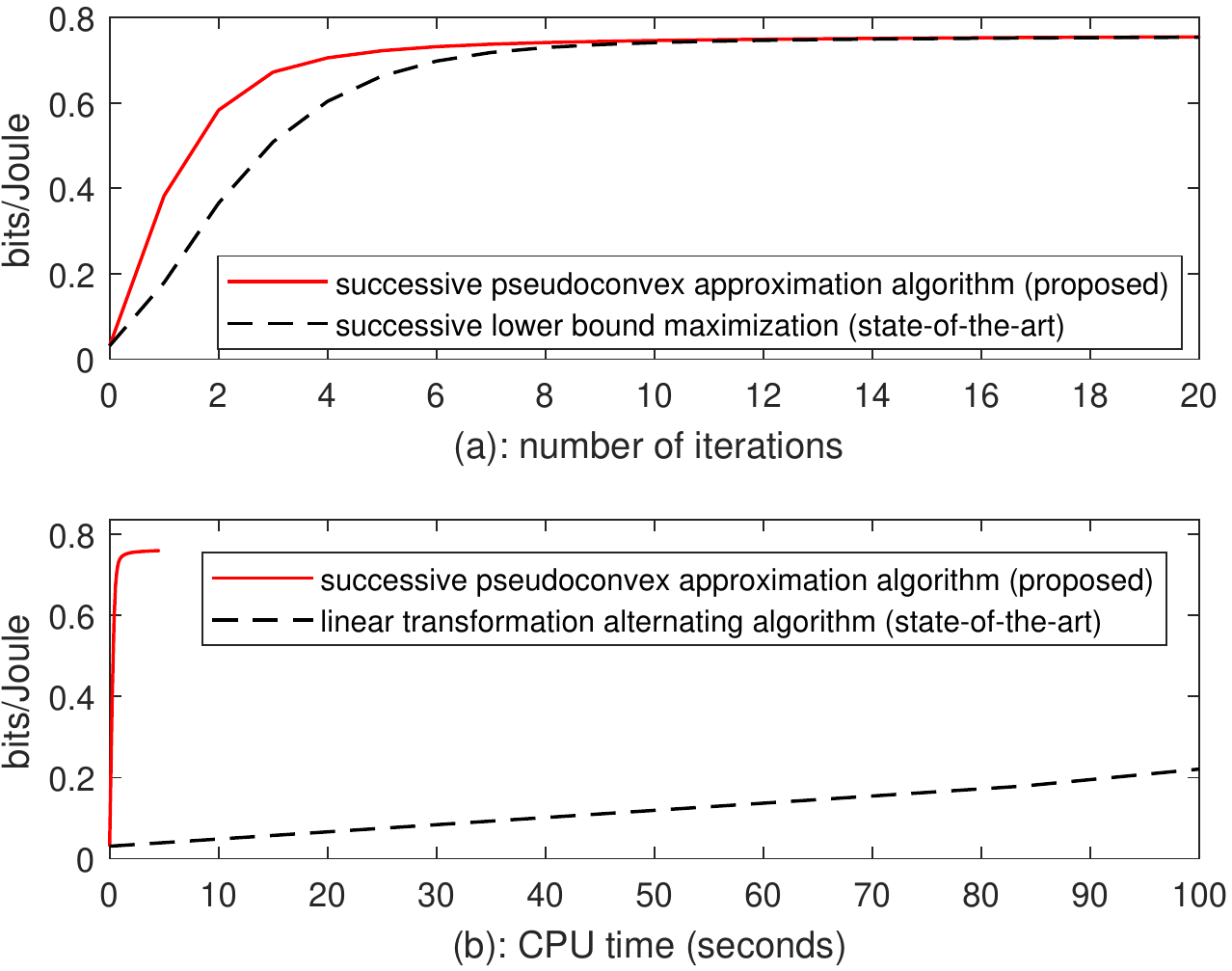}\caption{\label{fig:GEE-iteration}GEE maximization: achieved GEE vs. the number
of iterations and CPU time (seconds)}

\vspace{-2em}
\end{figure}

We compare the proposed Algorithm \ref{alg-GlobalEE} based on the
successive pseudoconvex approximation for problem (\ref{eq:original problem-GlobalEE})
with the successive lower bound minimization (SLBM) algorithm proposed
in \cite[Prop. 6]{Zappone2017}, which we briefly describe here. At
iteration $t$, $\mathbf{Q}^{t+1}$ is obtained by solving the following
problem:
\begin{equation}
\max_{\mathbf{Q}\succeq\mathbf{0}}\Biggl\{\negthinspace\frac{\sum_{k=1}^{K}(r_{k}^{+}(\mathbf{Q})-r_{k}^{-}(\mathbf{Q}^{t})-(\mathbf{Q}-\mathbf{Q}^{t})\bullet\nabla r_{k}^{-}(\mathbf{Q}^{t}))}{\sum_{k=1}^{K}P_{0,k}+\rho_{k}\textrm{tr}(\mathbf{Q}_{k})}\negthinspace\Biggr\},\label{eq:SLBM}
\end{equation}
subject to the power constraints ($\textrm{tr}(\mathbf{Q}_{k})\leq P_{k}$
for all $k$), where $r_{k}^{+}(\mathbf{Q})\triangleq\log\det(\sigma_{k}^{2}\mathbf{I}+\sum_{j=1}^{K}\mathbf{H}_{kj}\mathbf{Q}_{j}\mathbf{H}_{kj}^{H})$,
and this optimization problem is solved iteratively by the Dinkelbach's
algorithm. The SLBM algorithm bears its name from the fact that the
objective function in (\ref{eq:SLBM}) is a global lower bound of
the original objective function $f_{G}(\mathbf{Q})$ defined in (\ref{eq:original problem-GlobalEE}).
We do not consider the rate-dependent processing power consumption
here because the SLBM algorithm is not applicable otherwise.

As we see from Figure \ref{fig:GEE-iteration} (a), given the same
initial point ($\mathbf{Q}^{0}=P_{T}/M_{T}\mathbf{I}$), both algorithms
achieve the same GEE, and the proposed algorithm converges in fewer
number of iterations than the SLBM algorithm. However, as we see from
Figure \ref{fig:GEE-iteration} (b), the proposed algorithm needs
much less time to converge to a stationary point that that the SLBM
algorithm needs. This is because the variable update at each iteration
of the proposed algorithm can be implemented in closed-form expressions,
while a generic convex optimization problem in the form of (\ref{eq:SLBM})
must be solved (by CVX \cite{grant2011} in our simulations) for the
SLBM algorithm. Finally we remark that the SLBM algorithm cannot handle
rate-dependent processing power consumption.

\subsection{SEE Maximization}

We compare the proposed Algorithm \ref{alg-SumEE} based on the successive
pseudoconvex approximation for problem (\ref{eq:original-problem-SumEE})
with the linear transformation alternating (LTA) algorithm proposed
in \cite[Alg. 1]{He2014}. Note that we do not consider the rate-dependent
processing power consumption here because the LTA algorithm is not
applicable otherwise.

\begin{figure}[t]
\center\includegraphics[scale=0.6]{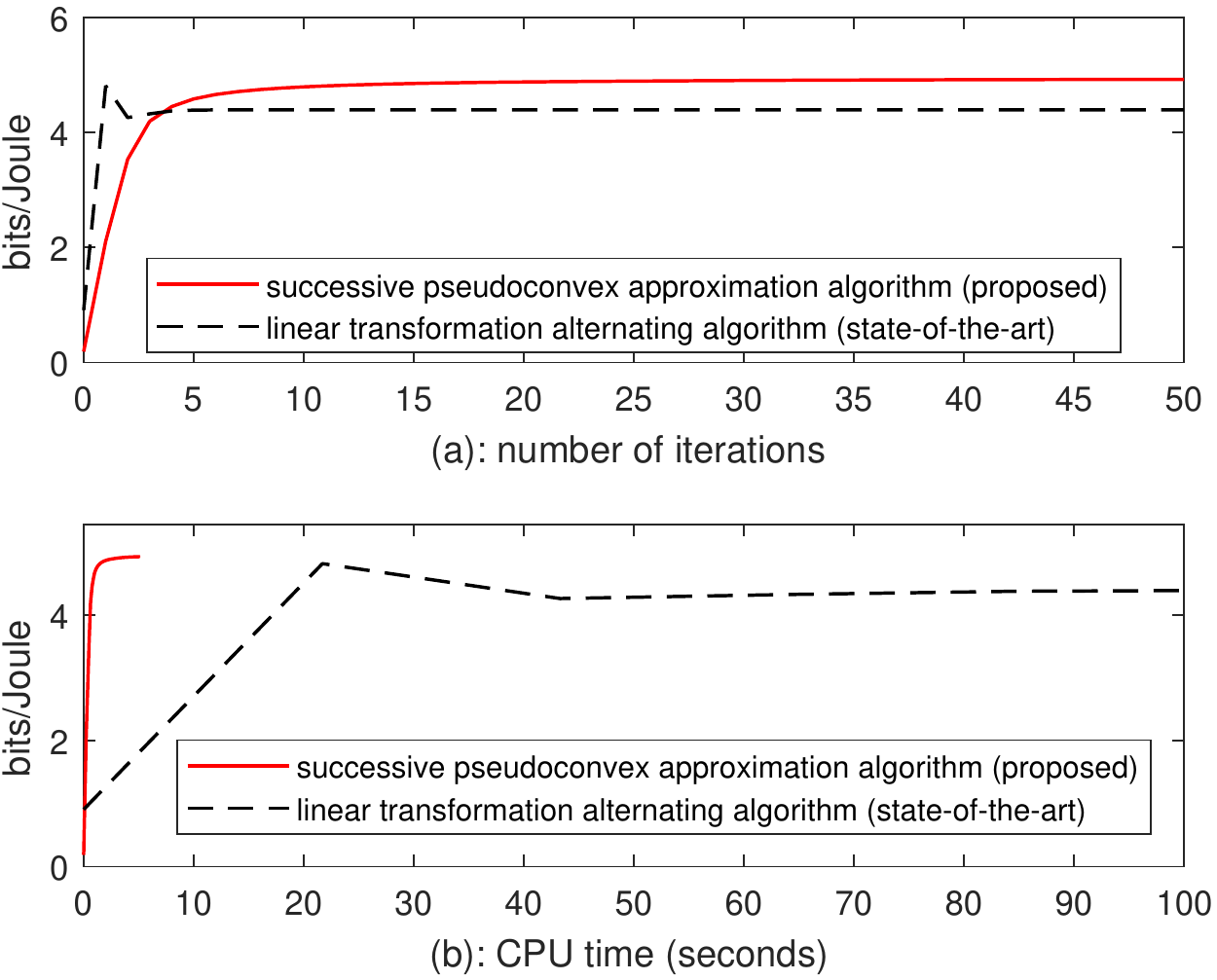}\caption{\label{fig:SEE-iteration}SEE maximization: achieved SEE vs. the number
of iterations and CPU time (seconds)}

\vspace{-2em}
\end{figure}

We can draw several observations from Figure \ref{fig:SEE-iteration},
where the achieved SEE versus the number of iterations and the CPU
time is plotted, respectively. Firstly, Figure \ref{fig:SEE-iteration}
(a) shows that the proposed algorithm achieves a better SEE than the
LTA algorithm. Secondly, as we can see from Figure \ref{fig:SEE-iteration}
(b), the proposed algorithm converges to the stationary point in less
than 1 second and is thus suitable for real time applications. This
is because the variable update at each iteration of the proposed algorithm
can be implemented in closed-form expressions and the Dinkelbach's
algorithm in the inner converges superlinearly. Although the variable
updates of the LTA algorithm are also based on closed-form expressions,
the inner layer is a BCD type algorithm which suffers from slow asymptotic
convergence and typically needs many iterations before convergence.
Finally we remark that the LTA algorithm cannot handle rate-dependent
processing power consumption.

\subsection{GEE and SEE maximization with QoS constraints}

\begin{figure}[t]
\center\includegraphics[scale=0.58]{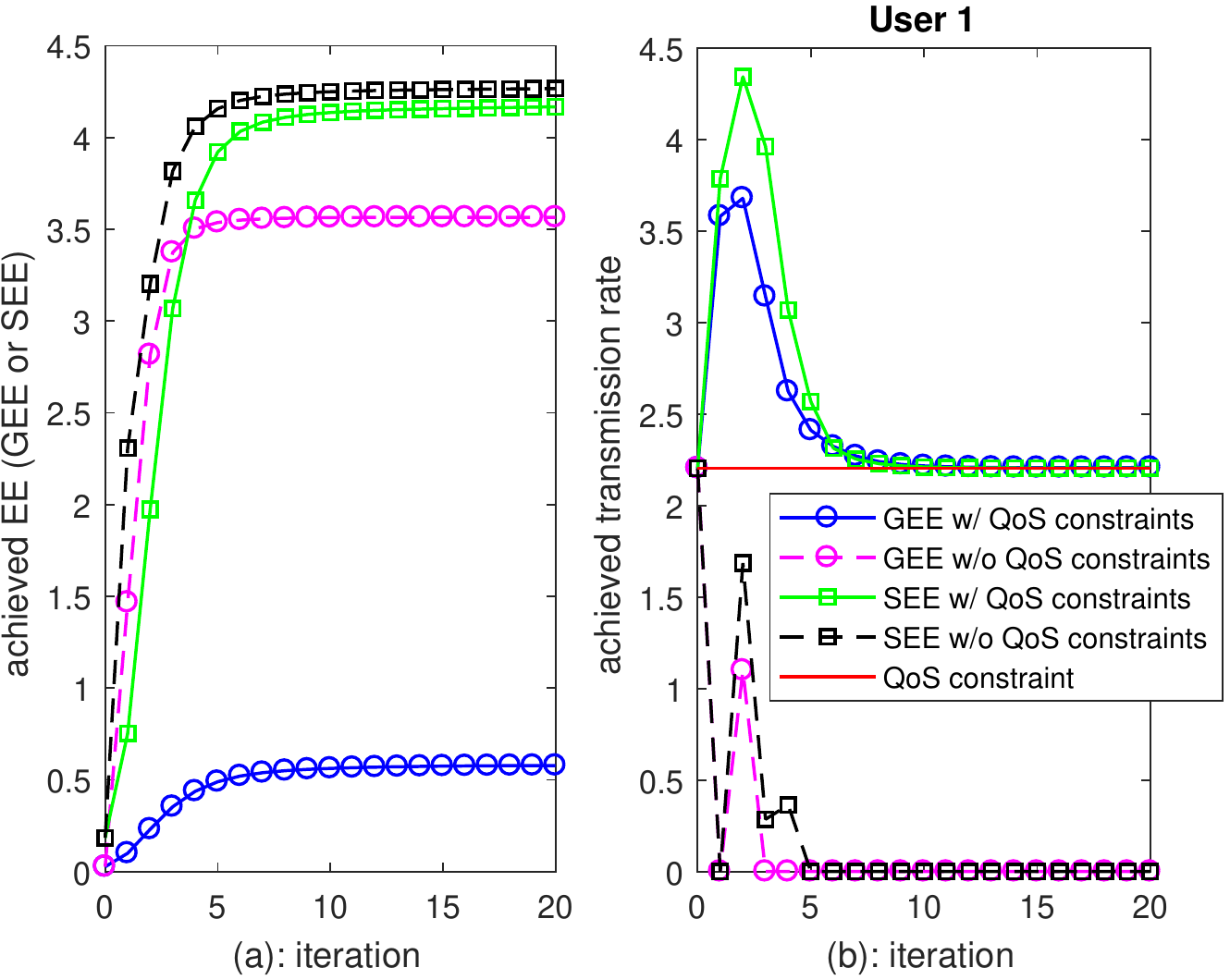}\caption{\label{fig:EE-wQoS}GEE and SEE maximization with QoS constraints:
achieved EE and achieved transmission rate of a particular user vs.
the number of iterations}

\vspace{-2em}
\end{figure}

In this subsection, we test Algorithm \ref{alg-GEE-w/QoS}-\ref{alg-SEE-w/QoS}
for problems (\ref{eq:original-problem-GlobalEE-w/QoS}) and (\ref{eq:original-problem-SumEE-w/QoS}).
In particular, the achieved GEE by Algorithm \ref{alg-GEE-w/QoS}
and the achieved SEE by Algorithm \ref{alg-SEE-w/QoS} is plotted
in Figure \ref{fig:EE-wQoS} (a). As a benchmark, we also plot the
achieved GEE by Algorithm \ref{alg-GlobalEE} and the achieved SEE
by Algorithm \ref{alg-SumEE}, which are designed for the GEE and
SEE maximization problems without QoS constraints, namely, (\ref{eq:original problem-GlobalEE})-(\ref{eq:original-problem-SumEE}).
On the one hand, we can see from Figure \ref{fig:EE-wQoS} (a) that
the achieved EE by Algorithms \ref{alg-GEE-w/QoS}-\ref{alg-SEE-w/QoS}
is, as expected, monotonically increasing w.r.t. the number of iterations.
The achieved GEE/SEE is smaller than that achieved by Algorithm \ref{alg-GlobalEE}/\ref{alg-SumEE},
because the feasible set of problem (\ref{eq:original-problem-GlobalEE-w/QoS})/(\ref{eq:original-problem-SumEE-w/QoS})
is only a subset of the feasible set of problem (\ref{eq:original problem-GlobalEE})/(\ref{eq:original-problem-SumEE}).
On the other hand, the transmission rate of a particular user is plotted
in Figure \ref{fig:EE-wQoS} (b) and we can see that this particular
user is guaranteed a minimum transmission rate by Algorithms \ref{alg-GEE-w/QoS}-\ref{alg-SEE-w/QoS},
while such a guarantee is not provided by Algorithms \ref{alg-GlobalEE}-\ref{alg-SumEE}.
This is because as long as the QoS constraints are not enforced, the
users with bad channel conditions may not be able to transmit in order
to maximize the GEE/SEE.

\section{\label{sec:Concluding-Remarks}Concluding Remarks}

In this paper, we have proposed novel iterative algorithms based on
the successive pseudoconvex approximation framework for the GEE and
SEE maximization problem, possibly with nonconcave QoS constraint
functions. As we have shown, pseudoconvexity plays a fundamental role,
because it enables us to design an approximate function that is not
necessarily a global lower bound of the original function. This makes
it possible to design new approximate functions that have more flexibility
(e.g., rate-dependent processing power consumption) and that can be
efficiently optimized. In particular, the proposed algorithms have
the following attractive features: 1) fast convergence as the structure
of the original optimization problem is preserved as much as possible
in the approximate problem solved in each iteration, 2) easy implementation
as each approximate problem is suitable for parallel computation and
its solution has a closed-form expression, and 3) guaranteed convergence
to a stationary point or a KKT point. These advantages of the proposed
algorithms are also numerically illustrated.

\appendix{}
\begin{IEEEproof}
[Proof of Proposition \ref{prop:ascent-direction-1}]Suppose $\mathbf{Q}^{t}=\mathbb{B}_{Q}\mathbf{Q}^{t}$.
The Lagrangian of (\ref{eq:approximate-problem-2}) is
\begin{align*}
 & \tilde{L}(\mathbf{Q},\mathbf{Y},\boldsymbol{\Pi},\boldsymbol{\Sigma},\boldsymbol{\lambda},\boldsymbol{\mu};\mathbf{Q}^{t},\mathbf{Y}^{t})\\
=\; & \tilde{f}_{G}(\mathbf{Q},\mathbf{Y};\mathbf{Q}^{t})+{\textstyle \sum_{k=1}^{K}}\boldsymbol{\Pi}_{k}\bullet\mathbf{Q}_{k}^{t}-{\textstyle \sum_{k=1}^{K}}\lambda_{k}(\textrm{tr}(\mathbf{Q}_{k})-P_{k})\\
 & +{\textstyle \sum_{k=1}^{K}}\mu_{k}(r_{k}^{+}(\mathbf{Y}_{k})-\overline{r}_{k}^{-}(\mathbf{Q}_{-k};\mathbf{Q}^{t}))\\
 & +{\textstyle \sum_{k=1}^{K}}\boldsymbol{\Sigma}_{k}\bullet(\mathbf{Y}_{k}-{\textstyle \sum_{j=1}^{K}}\mathbf{H}_{kj}\mathbf{Q}_{j}\mathbf{H}_{kj}^{H}),
\end{align*}
where $(\boldsymbol{\Pi},\boldsymbol{\Sigma},\boldsymbol{\mu},\boldsymbol{\lambda})$
are the dual variables. By definition $(\mathbf{Q}^{t},\mathbf{Y}^{t})$
solves the optimization problem (\ref{eq:approximate-problem-2}).
Since $(\mathbf{Q},\mathbf{Y})$ is a regular point \cite{Scutari_NonconvexApplications},
there exists $(\boldsymbol{\Pi}^{t},\boldsymbol{\Sigma}^{t},\boldsymbol{\mu}^{t},\boldsymbol{\lambda}^{t})$
such that $(\mathbf{Q}^{t},\mathbf{Y}^{t})$ and $(\boldsymbol{\Pi}^{t},\boldsymbol{\Sigma}^{t},\boldsymbol{\mu}^{t},\boldsymbol{\lambda}^{t})$
together satisfy the KKT conditions \cite[Prop. 4.3.1]{bertsekas1999nonlinear}:
\begin{equation}
\nabla_{\mathbf{Q}^{*}}\tilde{L}(\mathbf{Q}^{t},\mathbf{Y}^{t},\boldsymbol{\Pi}^{t},\boldsymbol{\Sigma}^{t},\boldsymbol{\lambda}^{t},\boldsymbol{\mu}^{t};\mathbf{Q}^{t},\mathbf{Y}^{t})=\mathbf{0},\label{eq:KKT-1}
\end{equation}
\begin{equation}
\mathbf{0}\preceq\boldsymbol{\Pi}_{k}^{t}\perp\mathbf{Q}_{k}^{t}\succeq\mathbf{0},0\leq\mu_{k}^{t}\perp\textrm{tr}(\mathbf{Q}_{k}^{t})-P_{k}\leq0,\label{eq:KKT-2}
\end{equation}
\begin{equation}
\mathbf{Y}_{k}^{t}={\textstyle \sum_{j=1}^{K}}\mathbf{H}_{kj}\mathbf{Q}_{j}^{t}\mathbf{H}_{kj}^{H},\boldsymbol{\Sigma}_{k}^{t}\bullet(\mathbf{Y}_{k}^{t}-{\textstyle \sum_{j=1}^{K}}\mathbf{H}_{kj}\mathbf{Q}_{j}^{t}\mathbf{H}_{kj}^{H}),\label{eq:KKT-3}
\end{equation}
\begin{equation}
0\leq\lambda_{k}^{t}\perp r_{k}^{+}(\mathbf{Y}_{k}^{t})-\overline{r}_{k}^{-}(\mathbf{Q}_{-k}^{t};\mathbf{Q}^{t})-R_{k}\geq0,\:\forall k.\label{eq:KKT-4}
\end{equation}
Substituting (\ref{eq:equal-gradient}) and (\ref{eq:lower-bound-1})
into (\ref{eq:KKT-1}) yields
\begin{equation}
\nabla_{\mathbf{Q}^{*}}L(\mathbf{Q}^{t},\mathbf{Y}^{t},\boldsymbol{\Pi}^{t},\boldsymbol{\Sigma}^{t},\boldsymbol{\lambda}^{t},\boldsymbol{\mu}^{t})=\mathbf{0},\label{eq:KKT-5}
\end{equation}
where $L(\mathbf{Q},\mathbf{Y},\boldsymbol{\Pi},\boldsymbol{\Sigma},\boldsymbol{\lambda},\boldsymbol{\mu})$
is the Lagrangian of (\ref{eq:GEE-QoS-new-problem}):
\begin{align*}
 & L(\mathbf{Q},\mathbf{Y},\boldsymbol{\Pi},\boldsymbol{\Sigma},\boldsymbol{\lambda},\boldsymbol{\mu})\\
=\; & f_{G}(\mathbf{Q})+{\textstyle \sum_{k=1}^{K}}\boldsymbol{\Pi}_{k}\bullet\mathbf{Q}_{k}^{t}-{\textstyle \sum_{k=1}^{K}}\lambda_{k}(\textrm{tr}(\mathbf{Q}_{k})-P_{k})\\
 & +{\textstyle \sum_{k=1}^{K}}\mu_{k}(r_{k}^{+}(\mathbf{Y}_{k})-r_{k}^{-}(\mathbf{Q}))\\
 & +{\textstyle \sum_{k=1}^{K}}\boldsymbol{\Sigma}_{k}\bullet(\mathbf{Y}_{k}-{\textstyle \sum_{j=1}^{K}}\mathbf{H}_{kj}\mathbf{Q}_{j}\mathbf{H}_{kj}^{H})
\end{align*}
Similarly, substituting (\ref{eq:lower-bound-1}) into (\ref{eq:KKT-4})
yields
\begin{equation}
0\leq\lambda_{k}\perp r_{k}^{+}(\mathbf{Y}_{k}^{t})-r_{k}^{-}(\mathbf{Q}^{t})-R_{k}\geq0,\:\forall k.\label{eq:KKT-6}
\end{equation}
Therefore $(\mathbf{Q}^{t},\mathbf{Y}^{t},\boldsymbol{\Pi}^{t},\boldsymbol{\Sigma}^{t},\boldsymbol{\lambda}^{t},\boldsymbol{\mu}^{t})$
satisfies the KKT conditions of problem (\ref{eq:GEE-QoS-new-problem}),
namely, (\ref{eq:KKT-2}), (\ref{eq:KKT-3}), (\ref{eq:KKT-5}) and
(\ref{eq:KKT-6}).

If, reversely, there exist $(\boldsymbol{\Pi}^{t},\boldsymbol{\mu}^{t},\boldsymbol{\lambda}^{t})$
and $(\mathbf{Q}^{t},\mathbf{Y}^{t})$ satisfying the KKT conditions
of problem (\ref{eq:GEE-QoS-new-problem}), namely, (\ref{eq:KKT-2}),
(\ref{eq:KKT-3}), (\ref{eq:KKT-5}) and (\ref{eq:KKT-6}), we can
see that $(\boldsymbol{\Pi}^{t},\boldsymbol{\mu}^{t},\boldsymbol{\lambda}^{t})$
and $(\mathbf{Q}^{t},\mathbf{Y}^{t})$ also satisfies the KKT conditions
of (\ref{eq:approximate-problem-2}), namely, (\ref{eq:KKT-1})-(\ref{eq:KKT-4}).
Since the objective function in (\ref{eq:approximate-problem-2})
is pseudoconcave and the constraint set $\bar{\mathcal{Q}}(\mathbf{Q}^{t})$
is convex, it follows from \cite[Th. 10.1.1]{Mangasarian_NonlinearProgramming}
that $(\mathbf{Q}^{t},\mathbf{Y}^{t})$ is an optimal point of (\ref{eq:approximate-problem-2}),
i.e., $\mathbb{B}_{Q}\mathbf{Q}^{t}=\mathbf{Q}^{t}$.

If $\mathbb{B}_{Q}\mathbf{Q}^{t}\neq\mathbf{Q}^{t}$, then
\[
\tilde{f}_{G}(\mathbb{B}_{Q}\mathbf{Q}^{t},\mathbb{B}_{Y}\mathbf{Q}^{t};\mathbf{Q}^{t})<\tilde{f}_{G}(\mathbf{Q}^{t},\mathbf{Y}^{t};\mathbf{Q}^{t}).
\]
Since $\tilde{f}_{G}(\mathbf{Q,Y};\mathbf{Q}^{t})$ is pseudoconcave,
\begin{align*}
0<\; & (\mathbb{B}_{Q}\mathbf{Q}^{t}-\mathbf{Q}^{t})\bullet\nabla_{\mathbf{Q}^{*}}\tilde{f}_{G}(\mathbf{Q}^{t},\mathbf{Y}^{t};\mathbf{Q}^{t})\\
 & +(\mathbb{B}_{Y}\mathbf{Q}^{t}-\mathbf{Y}^{t})\bullet\nabla_{\mathbf{Y}^{*}}\tilde{f}_{G}(\mathbf{Q}^{t},\mathbf{Y}^{t};\mathbf{Q}^{t})\\
= & (\mathbb{B}_{Q}\mathbf{Q}^{t}-\mathbf{Q}^{t})\bullet\nabla_{\mathbf{Q}^{*}}f_{G}(\mathbf{Q}^{t}),
\end{align*}
where the equality follows from (\ref{eq:equal-gradient}). Thus $\mathbb{B}_{Q}\mathbf{Q}^{t}-\mathbf{Q}^{t}$
is an ascent direction of $f_{G}(\mathbf{Q})$ at $\mathbf{Q}=\mathbf{Q}^{t}$.
\end{IEEEproof}


\begin{thebibliography}{10}
\providecommand{\url}[1]{#1}
\csname url@samestyle\endcsname
\providecommand{\newblock}{\relax}
\providecommand{\bibinfo}[2]{#2}
\providecommand{\BIBentrySTDinterwordspacing}{\spaceskip=0pt\relax}
\providecommand{\BIBentryALTinterwordstretchfactor}{4}
\providecommand{\BIBentryALTinterwordspacing}{\spaceskip=\fontdimen2\font plus
\BIBentryALTinterwordstretchfactor\fontdimen3\font minus
  \fontdimen4\font\relax}
\providecommand{\BIBforeignlanguage}[2]{{%
\expandafter\ifx\csname l@#1\endcsname\relax
\typeout{** WARNING: IEEEtran.bst: No hyphenation pattern has been}%
\typeout{** loaded for the language `#1'. Using the pattern for}%
\typeout{** the default language instead.}%
\else
\language=\csname l@#1\endcsname
\fi
#2}}
\providecommand{\BIBdecl}{\relax}
\BIBdecl

\bibitem{IMT}
ITU-R, ``{IMT vision -- Framework and oveall objectives of the future
  development of IMT for 2020 and beyond},'' Tech. Rep., 9 2015.

\bibitem{Luo2008}
\BIBentryALTinterwordspacing
Z.-Q. Luo and S.~Zhang, ``{Dynamic Spectrum Management: Complexity and
  Duality},'' \emph{IEEE Journal of Selected Topics in Signal Processing},
  vol.~2, no.~1, pp. 57--73, 2008.
\BIBentrySTDinterwordspacing

\bibitem{Ng2012}
D.~W.~K. Ng, E.~S. Lo, and R.~Schober, ``{Energy-efficient resource allocation
  in multi-cell OFDMA systems with limited backhaul capacity},'' \emph{IEEE
  Transactions on Wireless Communications}, vol.~11, no.~10, pp. 3618--3631,
  2012.

\bibitem{Ng2012a}
------, ``{Energy-efficient resource allocation in OFDMA systems with large
  numbers of base station antennas},'' \emph{IEEE Transactions on Wireless
  Communications}, vol.~11, no.~9, pp. 3292--3304, 2012.

\bibitem{Xu2013}
\BIBentryALTinterwordspacing
J.~Xu and L.~Qiu, ``{Energy Efficiency Optimization for MIMO Broadcast
  Channels},'' \emph{IEEE Transactions on Wireless Communications}, vol.~12,
  no.~2, pp. 690--701, feb 2013.
\BIBentrySTDinterwordspacing

\bibitem{Xu2014}
Q.~Xu, X.~Li, H.~Ji, and X.~Du, ``{Energy-efficient resource allocation for
  heterogeneous services in OFDMA downlink networks: Systematic perspective},''
  \emph{IEEE Transactions on Vehicular Technology}, vol.~63, no.~5, pp.
  2071--2082, 2014.

\bibitem{Tervo2015}
O.~Tervo, L.~N. Tran, and M.~Juntti, ``{Optimal Energy-Efficient Transmit
  Beamforming for Multi-User MISO Downlink},'' \emph{IEEE Transactions on
  Signal Processing}, vol.~63, no.~20, pp. 5574--5588, 2015.

\bibitem{Tervo2017}
\BIBentryALTinterwordspacing
O.~Tervo, A.~Tolli, M.~Juntti, and L.-N. Tran, ``{Energy-Efficient Beam
  Coordination Strategies With Rate-Dependent Processing Power},'' \emph{IEEE
  Transactions on Signal Processing}, vol.~65, no.~22, pp. 6097--6112, nov
  2017.
\BIBentrySTDinterwordspacing

\bibitem{Nguyen2015}
K.~G. Nguyen, L.~N. Tran, O.~Tervo, Q.~D. Vu, and M.~Juntti, ``{Achieving
  energy efficiency fairness in multicell MISO downlink},'' \emph{IEEE
  Communications Letters}, vol.~19, no.~8, pp. 1426--1429, 2015.

\bibitem{Pan2016}
C.~Pan, W.~Xu, J.~Wang, H.~Ren, W.~Zhang, N.~Huang, and M.~Chen,
  ``{Pricing-Based Distributed Energy-Efficient Beamforming for MISO
  Interference Channels},'' \emph{IEEE Journal on Selected Areas in
  Communications}, vol.~34, no.~4, pp. 710--722, 2016.

\bibitem{Zappone2017}
A.~Zappone, E.~Bj{\"{o}}rnson, L.~Sanguinetti, and E.~Jorswieck, ``{Globally
  Optimal Energy-Efficient Power Control and Receiver Design in Wireless
  Networks},'' \emph{IEEE Transactions on Signal Processing}, vol.~65, no.~11,
  pp. 2844--2859, 2017.

\bibitem{Rost2010}
P.~Rost and G.~Fettweis, ``{On the transmission-computation-energy tradeoff in
  wireless and fixed networks},'' \emph{2010 IEEE Globecom Workshops, GC'10},
  pp. 1394--1399, 2010.

\bibitem{Bjornson2015}
\BIBentryALTinterwordspacing
E.~Bj{\"{o}}rnson, L.~Sanguinetti, J.~Hoydis, and M.~Debbah, ``{Optimal Design
  of Energy-Efficient Multi-User MIMO Systems: Is Massive MIMO the Answer?}''
  \emph{IEEE Transactions on Wireless Communications}, vol.~14, no.~6, pp.
  3059--3075, jun 2015.
\BIBentrySTDinterwordspacing

\bibitem{He2014}
S.~He, Y.~Huang, L.~Yang, and B.~Ottersten, ``{Coordinated multicell multiuser
  precoding for maximizing weighted sum energy efficiency},'' \emph{IEEE
  Transactions on Signal Processing}, vol.~62, no.~3, pp. 741--751, 2014.

\bibitem{Yang_ConvexApprox}
Y.~Yang and M.~Pesavento, ``{A Unified Successive Pseudoconvex Approximation
  Framework},'' \emph{IEEE Transactions on Signal Processing}, vol.~65, no.~13,
  pp. 3313--3328, 2017.

\bibitem{Zappone2015}
\BIBentryALTinterwordspacing
A.~Zappone, L.~Sanguinetti, G.~Bacci, E.~Jorswieck, and M.~Debbah,
  ``{Energy-Efficient Power Control: A Look at 5G Wireless Technologies},''
  \emph{IEEE Transactions on Signal Processing}, vol.~64, no.~7, pp.
  1668--1683, apr 2016.
\BIBentrySTDinterwordspacing

\bibitem{Yang_MaxMinEE}
\BIBentryALTinterwordspacing
Y.~Yang and M.~Pesavento, ``{Energy efficiency in MIMO interference channels:
  Social optimality and max-min fairness},'' to appear in Proc. 2018
  International Conference On Acoustics, Speech and Signal Processing (ICASSP).
  [Online]. Available: \url{http://orbilu.uni.lu/handle/10993/34234}
\BIBentrySTDinterwordspacing

\bibitem{EarthD2.3}
\BIBentryALTinterwordspacing
``Energy efficiency analysis of the reference systems, areas of improvements
  and target breakdown,'' {INFSO-ICT-247733 EARTH D2.3}. [Online]. Available:
  \url{https://bscw.ict-earth.eu/pub/bscw.cgi/d71252/EARTH_WP2_D2.3_v2.pdf}
\BIBentrySTDinterwordspacing

\bibitem{Mangasarian_NonlinearProgramming}
O.~L. Mangasarian, \emph{{Nonlinear programming}}.\hskip 1em plus 0.5em minus
  0.4em\relax McGraw-Hill, 1969.

\bibitem{bertsekas1999nonlinear}
D.~P. Bertsekas, \emph{{Nonlinear programming}}.\hskip 1em plus 0.5em minus
  0.4em\relax Athena Scientific, 1999.

\bibitem{Scutarib}
\BIBentryALTinterwordspacing
G.~Scutari, F.~Facchinei, P.~Song, D.~P. Palomar, and J.-S. Pang,
  ``{Decomposition by Partial Linearization: Parallel Optimization of
  Multi-Agent Systems},'' \emph{IEEE Transactions on Signal Processing},
  vol.~62, no.~3, pp. 641--656, feb 2014.
\BIBentrySTDinterwordspacing

\bibitem{Ortega&Rheinboldt}
J.~M. Ortega and W.~C. Rheinboldt, \emph{{Iterative solution of nonlinear
  equations in several variables}}.\hskip 1em plus 0.5em minus 0.4em\relax
  Academic, New York, 1970.

\bibitem{Yang_stochastic}
\BIBentryALTinterwordspacing
Y.~Yang, G.~Scutari, D.~P. Palomar, and M.~Pesavento, ``{A Parallel
  Decomposition Method for Nonconvex Stochastic Multi-Agent Optimization
  Problems},'' \emph{IEEE Transactions on Signal Processing}, vol.~64, no.~11,
  pp. 2949--2964, jun 2016.
\BIBentrySTDinterwordspacing

\bibitem{Scutari2017}
\BIBentryALTinterwordspacing
G.~Scutari, F.~Facchinei, and L.~Lampariello, ``{Parallel and Distributed
  Methods for Constrained Nonconvex Optimization--Part I: Theory},'' \emph{IEEE
  Transactions on Signal Processing}, vol.~65, no.~8, pp. 1929--1944, apr 2017.
\BIBentrySTDinterwordspacing

\bibitem{Freund1994}
\BIBentryALTinterwordspacing
R.~W. Freund and F.~Jarre, ``{An interior-point method for fractional programs
  with convex constraints},'' \emph{Mathematical Programming}, vol.~67, no.
  1-3, pp. 407--440, oct 1994.
\BIBentrySTDinterwordspacing

\bibitem{Scutari_NonconvexApplications}
\BIBentryALTinterwordspacing
G.~Scutari, F.~Facchinei, L.~Lampariello, P.~Song, and S.~Sardellitti,
  ``{Parallel and Distributed Methods for Nonconvex Optimization--Part II:
  Applications},'' vol.~65, no.~8, pp. 1945--1960, 2016.
\BIBentrySTDinterwordspacing

\bibitem{grant2011}
\BIBentryALTinterwordspacing
M.~Grant and S.~Boyd, ``{CVX: Matlab Software for Disciplined Convex
  Programming, version 2.0 beta},'' http://cvxr.com/cvx, 2013.
\BIBentrySTDinterwordspacing

\end{thebibliography}
\end{document}